\documentclass[12pt]{article}
\usepackage{amsfonts,amssymb,amsmath}

 \def\medskipamount{12pt} \def\smallskipamount{6pt}

\newcounter{bitcount}

\newcommand{\bit}[1]{\addtocounter{bitcount}{1}\pagebreak[3]\bigskip\medskip\noindent
{\bf \large #1}\nopagebreak\setcounter{equation}{0}\medskip\newline\noindent}

\renewcommand{\theequation}{\thesubsection .\arabic{equation}}
\renewcommand{\thesubsection}{\arabic{bitcount}}

\newcommand{\re}[1]{\mbox{\ref{#1}}}

\catcode`\@=\active
\catcode`\@=11
\def\@eqnnum{\hbox to .01pt{}\rlap{\bf \hskip -\displaywidth(\theequation)}}
\catcode`\@=12

\begin{document}


\catcode`\@=\active
\catcode`\@=11
\def\@eqnnum{\hb@xt@.01\p@{}%
                      \rlap{\normalfont\normalcolor
                        \hskip -\displaywidth(\theequation)}}

\newcommand{\nc}{\newcommand}


\nc{\bs}[1]{ \addvspace{\medskipamount} \pagebreak[3]
\refstepcounter{equation}
\noindent {\bf (\theequation) #1.} \begin{em} \nopagebreak}

\nc{\es}{\end{em} \par \addvspace{\medskipamount} } 

\nc{\ess}{\end{em} \par}

\nc{\br}[1]{ \addvspace{\medskipamount} \pagebreak[3]
\refstepcounter{equation} 
\noindent {\bf (\theequation) #1} \nopagebreak}

\nc{\brs}[1]{ \pagebreak[3]
\refstepcounter{equation} 
\noindent {\bf (\theequation) #1.} \nopagebreak}

\nc{\er}{\par \addvspace{\medskipamount} }


\nc{\vars}[2]
{{\mathchoice{\mb{#1}}{\mb{#1}}{\mb{#2}}{\mb{#2}}}}
\nc{\A}{\mathbb A}
\nc{\C}{{\mathbb C}} 
\nc{\Cx}{{\mathbb C}^\times}
\nc{\Gm}{{\mathbb G}_m}
\renewcommand{\H}{\mathbb H}
\nc{\N}{\mathbb N}
\nc{\Q}{\mathbb Q}
\nc{\R}{\mathbb R}
\nc{\Z}{\mathbb Z}
\renewcommand{\P}{\mathbb P} 
\renewcommand{\O}{\mathcal O}


\nc{\oper}[1]{\mathop{\mathchoice{\mbox{\rm #1}}{\mbox{\rm #1}}
{\mbox{\rm \scriptsize #1}}{\mbox{\rm \tiny #1}}}\nolimits}
\nc{\ad}{\oper{ad}}
\nc{\Ad}{\oper{Ad}}
\nc{\Aut}{\oper{Aut}}
\nc{\aut}{\oper{aut}}
\nc{\Char}{\oper{char}}
\nc{\diag}{\oper{diag}}
\nc{\Ext}{\oper{Ext}}
\nc{\Hom}{\oper{Hom}}
\nc{\id}{\oper{id}}
\nc{\Lie}{\oper{Lie}}
\renewcommand{\mod}{\oper{mod}}
\nc{\Pic}{\oper{Pic}\,}
\nc{\Spec}{\oper{Spec}\,}
\nc{\Tor}{\oper{Tor}}

\nc{\operlim}[1]{\mathop{\mathchoice{\mbox{\rm #1}}{\mbox{\rm #1}}
{\mbox{\rm \scriptsize #1}}{\mbox{\rm \tiny #1}}}}

\nc{\lp}{\raisebox{-.1ex}{\rm\large(}}
\nc{\rp}{\raisebox{-.1ex}{\rm\large)}}


\nc{\al}{\alpha}
\nc{\be}{\beta}
\nc{\la}{\lambda}
\nc{\La}{\Lambda}
\nc{\ep}{\varepsilon}
\nc{\si}{\sigma}
\nc{\om}{\omega}
\nc{\Om}{\Omega}
\nc{\Ga}{\Gamma}
\nc{\Si}{\Sigma}


\nc{\Left}[1]{\hbox{$\left#1\vbox to
    11.5pt{}\right.\nulldelimiterspace=0pt \mathsurround=0pt$}}
\nc{\Right}[1]{\hbox{$\left.\vbox to
    11.5pt{}\right#1\nulldelimiterspace=0pt \mathsurround=0pt$}}
\nc{\updown}{\hbox{$\left\updownarrow\vbox to
    10pt{}\right.\nulldelimiterspace=0pt \mathsurround=0pt$}}


\nc{\beqas}{\begin{eqnarray*}}
\nc{\co}{{\cal O}}
\nc{\cx}{{\C^{\times}}}
\nc{\down}{\Big\downarrow}
\nc{\Down}{\left\downarrow
    \rule{0em}{8.5ex}\right.}
\nc{\downarg}[1]{{\phantom{\scriptstyle #1}\Big\downarrow
    \raisebox{.4ex}{$\scriptstyle #1$}}}
\nc{\leftdownarg}[1]{\raisebox{.4ex}{$\scriptstyle #1$}\Big\downarrow
    {\phantom{\scriptstyle #1}}}
\nc{\eeqas}{\end{eqnarray*}}
\nc{\fp}{\mbox{     $\Box$} \par \addvspace{\smallskipamount}}
\nc{\lrow}{\longrightarrow}
\nc{\pf}{\noindent {\em Proof}}
\nc{\sans}{\, \backslash \,}
\nc{\st}{\, | \,}


\nc{\Exist}[1]{Extension of structure group induces a sur\-jec\-tion
  $H^1(#1,T) \to \bar{H}^1(#1,G)$; that is, every rationally trivial
  principal $G$-bundle over $#1$ admits a reduction to $T$.}

\nc{\Unique}[1]{The latter surjection is also an injection: that is,
  for every rationally trivial principal $G$-bundle $E$ over $#1$, the
  isomorphism class of the reduction to $T$ is unique modulo the
  action of $W$.}

\nc{\Corollary}[1]{Let $P = \Pic #1$; then the surjection above
  descends to a natural surjection $(P \otimes \Lambda)/W \to
  \bar{H}^1(#1,G)$.}

\nc{\Connect}[1]{automorphism group $\Ga (#1, \Ad E)$ is smooth,
  affine, and connected, as is the kernel of the evaluation map $\Ga
  (#1, \Ad E) \to G^\ell$ defined when $E$ is trivialized at $\ell$
  rational points.}

\hyphenation{para-met-riz-ing sub-bundle sur-jec-tion cen-tral-iz-er normal-iz-er}

\catcode`\@=12



\noindent
{\LARGE \bf Variations on a theme of Grothendieck} \smallskip \\
{\large Or, I've a feeling we're not in Kansas any more}\\ 

\noindent
{\bf Johan Martens } \\ 
School of Mathematics and Maxwell Institute, University of Edinburgh\\
\verb+johan.martens@ed.ac.uk+

\smallskip

\noindent
{\bf Michael Thaddeus }  \\ 
Department of Mathematics, Columbia University \\
\verb+thaddeus@math.columbia.edu+

\renewcommand{\thefootnote}{} 

\footnotetext{JM was supported by QGM (Centre for Quantum Geometry of
  Moduli Spaces), funded by the Danish National Research Foundation,
  and was partially supported by the ESF Research Networking Programme
  ``Low-Dimensional Topology and Geometry with Mathematical Physics''
  (ITGP).}

\footnotetext{MT was partially supported by NSF grants DMS--0401128
  and DMS--0700419 and wishes to thank the University of Geneva and
  Imperial College, London for their hospitality while this paper was
  written.}

\setcounter{bitcount}{-1}

\bigskip

\noindent
The theme of the paper is that principal bundles over the complex
projective line with reductive structure group can be reduced to a
maximal torus, in a unique fashion modulo automorphisms and the action
of the Weyl group.  Grothendieck proved this in a remarkable paper
\cite{kansas} written, in quite a different style from that of his
later works, during the year 1955, which he spent in Kansas.
Harder \cite{harder} showed in 1968 that the statement remains valid
over an arbitrary ground field, provided that the reductive group is
split and the bundle is trivial over the generic point.

Our variations are extensions of this result, on the
existence and uniqueness of such a reduction, to principal bundles
over other base schemes or stacks:

\begin{enumerate}

\item a {\it $\mu_n$-equivariant line}, that is, a line modulo the
  group $\mu_n$ of $n$th roots of unity;

\item a {\it football}, that is, an orbifold whose coarse
  moduli space is the projective line, with orbifold structure
  at two points;

\item a {\it gerbe} with structure group $\mu_n$ over a football;

\item a {\it chain of lines} meeting in nodal singularities; 

\item a {\it torus-equivariant line}, that is, a line modulo an action of
  a split torus; and

\item a {\it torus-equivariant chain}, that is, a chain of lines
  modulo an action of a split torus.

\end{enumerate}

We also prove that the automorphism group schemes of such
bundles are smooth, affine, and connected.  All of the above
results, in the case of a torus-equivariant chain, are
essential to the authors' work on compactifications of
reductive groups \cite{wonderful}.

We are very grateful to Jarod Alper, Jack Evans, Robert Friedman, Jens
Carsten Jantzen, Martin Olsson, Zinovy Reichstein, Angelo Vistoli, and
the late Torsten Ekedahl for their invaluable help and advice.  Above
all, we wish to thank Brian Conrad, who was extremely helpful on
several occasions and who provided the proof of a useful result, on
the connectedness of centralizers in positive characteristic, which is
recounted in the Appendix.

\bit{Theme: the projective line}
Let $k$ be an arbitrary field.  All schemes and stacks throughout are over $k$.

Let $X$ be a scheme or stack, and let $G$, $H$ be linear algebraic
groups, that is, smooth affine group schemes of finite type over $k$.

A {\it principal bundle\/} over $X$ with {\it structure group} $G$ is
a $G$-scheme or $G$-stack $E \to X$ locally trivial in the \'etale
topology.  Isomorphism classes of principal $G$-bundles over $X$ are
naturally in bijection with the \'etale cohomology set $H^1(X,G)$.  

If $\rho: G \to H$ is a homomorphism, the {\it extension of structure
  group} is the principal $H$-bundle $E_\rho := [(E \times H)/G]$.  A
{\it reduction of structure group} to a subgroup $H \subset G$ is a
section of $[E/H]$, the fiber bundle associated to $E$ with fiber
$G/H$.  The inverse image of this section in $E$ is a principal
$H$-bundle, which we denote $E_H$; its extension of structure group by
the inclusion $H \subset G$ is naturally isomorphic to $E$.  Two
reductions of structure group of $E$ to $H$ may, of course, be interchanged
by an automorphism of the $G$-bundle $E$; but, as is easily seen, this
is the case if and only if the principal $H$-bundles are isomorphic.

Suppose henceforth that $G$ is a split reductive group, which we take to be
connected; let $T \subset G$ be a split maximal torus; and let $B
\subset G$ be a Borel subgroup containing $T$.  Let $\eta = \Spec
k(t)$ be the generic point of $\P^1$, and define a $G$-bundle over
$\P^1$ to be {\it rationally trivial} if its restriction to $\eta$ is
trivial.  The set, denoted $\bar{H}^1(\P^1,G)$, of isomorphism classes
of rationally trivial $G$-bundles is therefore nothing but $\ker
H^1(\P^1,G) \to H^1(\eta,G)$.  All $T$-bundles over $\P^1$ are
rationally trivial by Theorem 90 \cite[III 4.9]{milne}, so
$\bar{H}^1(\P^1,T) = H^1(\P^1,T)$.  Furthermore, if $k$ is
algebraically closed of characteristic zero, a theorem of Steinberg
\cite[1.9]{steinberg} implies that all $G$-bundles over $\P^1$ are
rationally trivial, so $\bar{H}^1(\P^1,G) = H^1(\P^1,G)$ in this case.

\bs{Theorem (Existence for lines)} 
\label{line-exist}
\Exist{\P^1}
\es

\noindent {\it Proof} following Harder \cite{harder}.

{\bf Step A}: {\it The case $G=PGL_2$}.  Let $F$ be a rationally
trivial $PGL_2$-bundle over $\P^1$.  The associated $\P^1$-bundle $[(F
\times \P^1)/PGL_2]$ is trivial over the generic point.  Let $D$ be
the closure of a generic section; then $\O(D)$ is a line bundle
restricting to $\O(1)$ on each fiber.  Its direct image under the
projection to $\P^1$ is then a rank 2 vector bundle whose dual $E$
satisfies $F = \P E$.  It therefore suffices to show that every rank 2
vector bundle $E$ over $\P^1$ is isomorphic to $\O(d) \oplus \O(e)$
for some $d \geq e$.

Any such bundle has a line subbundle $\O(d)$ of maximal degree.
Indeed, being rationally trivial, $E$ has a nonzero rational section.
Let $D$ be its divisor of zeroes and poles, and let $d = |D|$; there
is then a nowhere vanishing section of $E(-d)$ yielding a short exact
sequence of the form
$$0 \lrow \O(d) \lrow E \lrow \O(e) \lrow 0.$$ 
So a line subbundle exists.  Furthermore, if $\O(c) \to E$ is any line
subbundle, then either $c = d$ or there is a nonzero map $\O(c) \to
\O(e)$, in which case $c \leq e$.  Hence the set of possible $c$ is
bounded above; assume without loss of generality that $d$ is maximal
among them.

We claim that then $d \geq e$.  If not, then $H^0(\O(e-d-1)) \neq 0$,
so there exists a nonzero map $\O(d+1) \to \O(e)$.  Since in the long
exact sequence
$$H^0(\Hom(\O(d+1),E)) \lrow
H^0(\Hom(\O(d+1),O(e))) \lrow H^1(\Hom(\O(d+1),\O(d)))$$
the right-hand term vanishes, our nonzero map comes from a
nonzero map $\O(d+1) \to E$.  Let $D$ be its divisor of
vanishing; there is then a nowhere zero map $\O(d+1+|D|) \to
E$, contradicting the maximality of $d$.  This proves the
claim.

The extension class of the short exact sequence lies in
$H^1(\O(d-e))$, but by the claim this vanishes.  The sequence
therefore splits.  This completes Step A.

{\bf Step B}: {\it Reduction to $B$.}
Rational triviality guarantees that a generic section of the
associated $G/B$-bundle $[E/B] \to \P^1$ exists.  Any cover of $\P^1$
by finitely many \'etale trivializations of $E$ provides a
faithfully flat and quasi-compact morphism to $\P^1$ by which the base
change of $[E/B] \to \P^1$ is proper.  Hence $[E/B] \to \P^1$ is
itself proper \cite[2.7.1]{egaiv}, so the generic section extends to
a regular section by the valuative criterion.  Hence $E$ reduces to
$B$.

{\bf Step C}: {\it An adapted reduction to $B$.}  
We wish the corresponding $B$-bundle $E_B$ to enjoy the
property that for all positive roots $\alpha: B \to \Gm$, the
associated line bundle $E_{\alpha}$ has degree $\geq 0$.
This is false in general, but the given reduction to $B$ can be
modified to give a new one where this property holds.  We call
such a reduction an {\it adapted} reduction (Harder calls it ``reduziert'').

To prove this, observe first that it suffices for this property
to hold for the simple roots $\alpha_i$, since any positive
root is a nonnegative combination of them.  Next consider the
parabolic subgroup $P_i$ whose Lie algebra ${\mathfrak p}_i$ is
the sum of ${\mathfrak b}$ and the root space ${\mathfrak
  g}_{-\alpha_i}$.  Then $P_i/B \cong \P^1$, so the
$P_i/B$-bundle associated to $E_B$ is a $\P^1$-bundle with
structure group $PGL_2$.  Its total space is naturally
contained in that of the $G/B$-bundle above and contains the
chosen section.  But it also contains the section corresponding
to $\O(d)$, the line bundle of greater degree, in the splitting
of Step A\@.  If this section is different, it gives a
different reduction of $E$ to $B$.  For this new reduction, the
line bundle associated to $\alpha_i$ is $\O(d-e)$, whose degree
is nonnegative.  But the line bundles associated to all other
simple roots are unchanged from the previous reduction, since
for $j \neq i$ the homomorphisms $\alpha_j: B \to \Gm$
factor through $P_i$, and the section of the associated
$G/P_i$-bundle is unchanged.  We may therefore modify the
section of our $G/B$-bundle one simple root at a time until the
desired property holds for all simple roots.

{\bf Step D}: {\it Conclusion}.  Thus we obtain an adapted reduction
$E_B$ to $B$.  We seek a further reduction to $T$, that is, a section
of the associated $B/T$-bundle.  Recall now that the Borel is a
semidirect product $B = T \ltimes U$, where $U$ is the maximal
unipotent subgroup.  Hence $B/T$ is canonically isomorphic to $U$ as a
variety; but the natural actions of $B$ do not match, so the $B/T$-bundle
and the $U$-bundle associated to our $B$-bundle are not the same.
Rather, our $B/T$-bundle is a torsor for $F$, the bundle of groups
associated to $E_B$ via the conjugation action of $B$ on $U$.  The
identity of course furnishes a section of $F$, so it suffices to show
that all torsors for $F$ are trivial.

This is where we use the adapted property established in Step C\@.
As a group with $B$-action, $U$ is filtered by subgroups whose
successive quotients are isomorphic to the additive group ${\mathbb G}_a$,
with $B$-action given by a positive root $\alpha$.  Hence 
$F$ is filtered by subbundles of groups whose
successive quotients are line bundles of nonnegative degree.
On $\P^1$, any such line bundle has trivial $H^1$, so from the
long exact sequences associated to the filtration it follows
that $H^1(F) = 0$ and hence that every torsor over $F$ is
trivial.
\fp

\bs{Corollary}
\label{line-corollary}
Let $\La := \Hom(\Gm, T)$ be the cocharacter lattice, $W$ the Weyl
group of $G$.  Then there is a natural surjection from $\La/W$ to
$\bar{H}^1(\P^1,G)$.
\es

\pf.  Tensoring the usual isomorphism $\Z \to \Pic \P^1 =
H^1(\P^1,\Gm)$ by the free $\Z$-module $\La = \Hom(\Gm, T)$ yields an
isomorphism $\La \to H^1(\P^1,T)$.  By Theorem \re{line-exist}, its
composition with the extension of structure group $T \to G$ is a
surjection $\La \to H^1(\P^1,G)$.  There is some redundancy, however.
Since for any inner automorphism $\rho: G \to G$ there is a natural
isomorphism of principal bundles $E_\rho \cong E$, and every Weyl
group element $w:T \to T$ extends to an inner automorphism of $G$, it
follows that $\la$ and $w \cdot \la$ determine isomorphic $G$-bundles.
The surjection above therefore descends to $\La/W$.  \fp

\bs{Theorem (Uniqueness for lines)}
\label{line-unique}
\Unique{\P^1}
\es

\pf\/ following Grothendieck \cite{kansas}.  

First, we prove the uniqueness for the case $G = GL_n$, that is, for
vector bundles.  This amounts to proving that if $\bigoplus \O(a_i)
\cong \bigoplus \O(b_j)$, then the $a_i$ and the $b_j$ are the same up
to permutation.  Or equivalently, if both sequences are weakly
increasing, then they are the same.  Since only bundles of nonnegative
degree have nonzero sections, in order to have maps in both directions
nonvanishing on all summands, we must have $a_n = b_n$.  Then the map
from $\O(a_n)$ to $\bigoplus \O(b_j)$ must be a constant map to the
summands of degree $a_n$; we may split off its image and proceed by
induction on $n$ to finish the $GL_n$ case.

Next, we reduce the general case to the $GL_n$ case.  Suppose
$\la, \la' \in \La$ induce isomorphic $G$-bundles on $\P^1$.  
Then the vector bundles associated to any representation
$G \to GL_n$ are also isomorphic and hence split into line
bundles of the same degrees.  Conjugate the representation so that it takes
$T$ to the diagonal matrices.  Then its value on $\la(t)$ is
$\diag(t^{a_1} , \cdots , t^{a_n})$ where the integers $a_i$
are the splitting type, as above.  Hence the character of any 
representation takes the same values on $\la$ as on $\la'$.

By passing to a field extension, it suffices to assume that $k$ is
algebraically closed, with nonzero elements that are not roots of
unity.  Then the algebra of regular class functions on $G$ is
generated by the characters of irreducible representations
\cite[6.1]{steinberg}.  Hence all such functions agree on $\la$ and
$\la'$, so the compositions of $\la$ and $\la'$ with the projection to
$G/G = T/W$ are equal.  In particular, if $t_0 \in \Gm$ is a
$k$-rational point that is not a root of unity, so that it generates a
Zariski dense subset, then $\la(t_0)$ and $\la'(t_0)$ are conjugate by
some constant $w \in W$.  Hence the same is true of $\la$ and $\la'$
themselves.  \fp

Let us mention two further results, though they do not feature in
Grothendieck's paper: on the existence and uniqueness of a
Harder-Narasimhan reduction, and on the connectedness of the
automorphism group.

\bs{Theorem} 
\label{harder-narasimhan}
A rationally trivial $G$-bundle $E$ over $\P^1$ admits an unique
reduction to a parabolic subgroup $P$ containing $T$ which (a) is
rigid in that its infinitesimal deformations are trivial and (b) is
minimal among all such rigid reductions to parabolic subgroups.  Any
reduction of $E$ to $T$ is a further reduction of this one.  
\es

\pf.  Observe first that, as any parabolic $P$ is the normalizer in $G$ of
its own Lie algebra ${\mathfrak p} \subset {\mathfrak g}$
\cite[11.16]{borel}, $G/P$ parametrizes the subalgebras of ${\mathfrak
  g}$ conjugate to ${\mathfrak p}$.  Consequently, a reduction of $E$
to $P$, being a section $\si: \P^1 \to E/P$, determines, and is
determined by, a bundle $\ad E_P$ of subalgebras of $\ad E$ conjugate
to ${\mathfrak p}$.  The infinitesimal deformations of the reduction
are sections of the normal bundle to $\si(\P^1)$ in $E/P$, namely $\ad
E/\ad E_P$.

By Theorem \re{line-exist} there exists a reduction of $E$ to $T$.  Then by
Corollary \re{line-corollary}, $E$ is induced by a cocharacter $\la$,
so there is a splitting
$$\ad E \cong \O^r \oplus \bigoplus_{\al \in \Phi} \O(\al \cdot \la)$$
corresponding to the splitting of ${\mathfrak g}$ into Cartan
subalgebra and root spaces.  If any of the summands of nonnegative
degree is not contained in $\ad E_P$, then it will map nontrivially to
$\ad E/\ad E_P$, giving the latter a nonzero section.  Hence for the
reduction to be rigid, $\ad E_P$ must contain all of the summands of
nonnegative degree in the splitting above.

On the other hand, the summands of nonnegative degree span a bundle of
parabolic subalgebras with the desired rigidity.  It therefore
determines the unique minimal rigid reduction.  As $\ad E_P$ contains
the bundle of Cartan subalgebras $\ad E_T \cong \O^r$, this reduction
reduces further to $T$.  \fp

\bs{Theorem (Connectedness for lines)}
\label{line-connected} 
The \Connect{\P^1}
\es

\pf.  Let $E_P \subset E$ be the Harder-Narasimhan reduction of
Theorem \re{harder-narasimhan}.  Since it is unique and rigid, any
family of automorphisms must preserve it, so $\Ga(\Ad E) = \Ga (\Ad
E_P)$.  Let $L$ be the Levi factor of $P$ containing $T$ and let $U$
be the maximal unipotent, so that $P = L \ltimes U$.  The adjoint
action of $T$ preserves $L$ and $U$, and even the root subgroups
$U_\al$ that directly span $U$.  Hence $\Ad E_P$ becomes a semidirect
product of group bundles $Q \ltimes R$, where $Q$ is an adjoint
$L$-bundle, $R$ is an adjoint $U$-bundle, and $R$ is directly spanned
by line bundles $R_\al$.  Hence as a scheme $\Ga(\Ad E_P) \cong \Ga(Q)
\times \Ga(R)$, and furthermore $\Ga(R) \cong \prod_\al
\Ga(R_\al)$, where $\al$ runs over the root spaces in the Lie algebra of
$U$.

Consider first $\Ga(Q)$.  Let $\la: \Gm \to T$ be the 1-parameter
subgroup defining $E$, as usual, and note that for any root $\al$ of
the Levi, $\la \cdot \al = 0$.  Hence $\la(\Gm)$ is contained in the
center of the Levi, so its adjoint action on $L$ is trivial, and $Q$
is trivial so that $\Ga(Q) = L$, which is smooth, affine, and
connected.

On the other hand, $\Ga(R) \cong \prod_\al \Ga(R_\al)$, and since
$R_\al$ is a line bundle, $\Ga(R_\al)$ is an affine space, so $\Ga(R)$
is also an affine space, hence certainly smooth, affine, and connected.

As for the kernel of the evaluation map, it is a semidirect product of
the corresponding kernels for $\Ga(Q)$ and $\Ga(R)$.  From the above,
the first kernel is clearly trivial, while the second one is
clearly an affine space.  \fp

\bit{\boldmath Variation \thesubsection: a $\mu_n$-equivariant line}
Let $\mu_n$ be the group of $n$th roots of unity, acting on the
projective line $\P^1$ by $\la [x,y] = [\la x,y]$.  This
variation is concerned with $G$-bundles over the stack $[\P^1/\mu_n]$,
or equivalently, with $\mu_n$-equivariant $G$-bundles over $\P^1$.

A $G$-bundle over $[\P^1/\mu_n]$ is said to be {\it rationally
  trivial\/} if its restriction to the generic point $\Spec k(t)$ of
$[\P^1/\mu_n]$ is trivial.

\br{Counterexample}\label{mu-line-example}where the bundle is not rationally
trivial and does not reduce to $T$.  If $k$ is the field of real
numbers and $G = PGL_2$, then the subgroup generated by
$\left( \begin{smallmatrix} \phantom{-}0&1\\-1&0 \end{smallmatrix}
\right)$ is isomorphic to $\mu_2$ but lies in no split torus.  For all
split tori are conjugate over $k$ \cite[20.9]{borel}, so the
eigenvalues of any matrix in a split torus are real.  The bundle
$[(\P^1 \times G)/\mu_2]$ over $[\P^1/\mu_2]$ therefore does not
reduce to the split torus $T$.  In light of the following lemma, this
bundle cannot be rationally trivial. \er

\bs{Lemma}
\label{mu-line-lemma}
A homomorphism $\phi: \mu_n \to PGL_2$ has image in a split torus if
and only if the induced $PGL_2$-bundle $[(\P^1 \times PGL_2)/\mu_n]$
over $[\P^1/\mu_n]$ is rationally trivial.
\es

\pf. For such a bundle to be rationally trivial means that there
exists a rational map $f: \Gm \dashrightarrow PGL_2$ such that the
following diagram commutes, where the top arrow is multiplication and the
bottom arrow is given by $(\la,g) \mapsto \phi(\la)g$:
$$
\begin{array}{ccc}
\mu_n \times \Gm & \lrow & \Gm \\
\downarg{\id \times f} && \downarg{f} \\
\mu_n \times PGL_2 & \lrow & PGL_2.
\end{array}
$$
If $\phi$ has image in a split torus, then it extends to a 1-parameter
subgroup $f: \Gm \to PGL_2$, which has the desired property.

Conversely, suppose that the rational map $f$ exists.  Consider the
completion of $PGL_2$ to the projective space $\P^3$ parametrizing all
$2 \times 2$ matrices modulo scalars.  By the valuative criterion, $f$
extends to a morphism $\bar{f}: \P^1 \to \P^3$, and the diagram
$$
\begin{array}{ccc}
\mu_n \times \P^1 & \lrow & \P^1 \\
\downarg{\id \times \bar{f}} && \downarg{\bar{f}} \\
\mu_n \times \P^3 & \lrow & \P^3
\end{array}
$$
still commutes.

Every torus in $PGL_2$ fixes two points in $\P^1$ over some
extension field of $k$.  If one of these points is $k$-rational, then
clearly so is the other.  As the split tori are all conjugate over $k$
\cite[20.9]{borel}, this occurs if and only if the torus is split.
Hence $\phi$ has image in a split torus if and only if it fixes
some $k$-rational point in $\P^1$.

If $\phi$ is trivial, the statement is immediate.  Otherwise, the
action of $\mu_n$ on $PGL_2$ fixes no points, while that on the
quadric surface $\P^3 \sans PGL_2 \cong \P^1 \times \P^1$ is induced
from the action of $\mu_n$ on the first factor via $\phi$.  Since the
action of $\mu_n$ on $\P^1$ by $\la[x,y] = [\la x,y]$ fixes $[1,0]$
and $[0,1]$, the commutativity of the diagram implies they must map to
$k$-rational points on the quadric whose first factors are fixed by
the $\mu_n$-action via $\phi$.  Hence $\phi$ has image in a split
torus.  \fp

\bs{Theorem (Existence for {\boldmath $\mu_n$}-equivariant lines)}
\label{mu-line-exist}
\Exist{[\P^1/\mu_n]}
\es

\pf. 
{\bf Step A}: {\it The case $G =PGL_2$.}  Let $E$ be a rationally trivial
$PGL_2$-bundle over $[\P^1/\mu_n]$.  By Theorem \re{line-exist}, its
pullback to $\P^1$ lifts to a $GL_2$-bundle which splits as $\O(d)
\oplus \O(e)$ with $d \geq e$.  If $d = e$, then $E$ is trivial as a
$PGL_2$-bundle over $\P^1$; since any split torus is conjugate to $T$,
the result then follows from Lemma \re{mu-line-lemma}.  

Otherwise, the section of $\pi: \P (\O(d) \oplus \O(e)) \to \P^1$
corresponding to $\O(d)$ is the unique section with negative normal
bundle and hence is preserved by the $\mu_n$-action.  It is therefore
a $\mu_n$-invariant effective divisor $D$ restricting to $\O(1)$ on
each fiber of $\pi$.  Hence $\O(D)$ is an equivariant line bundle,
that is, a line bundle over $[\P^1/\mu_n]$, so that $\tilde{E} =
\pi_*\O(D)^*$ is an equivariant lifting of $E$ to $GL_2$ preserving
$\O(1)$.  That is, over $[\P^1/\mu_n]$ there is an exact sequence 
$$0 \lrow \O(d) \lrow \tilde{E} \lrow \O(e) \lrow 0.$$ 

The extension class lies in $H^1([\P^1/\mu_n],\O(d-e))$.  Since
$\mu_n$ is a multiplicative group, its rational representations are
completely reducible, so taking invariants is exact.  Hence any
line bundle $L$ over $[\P^1/\mu_n]$ satisfies $H^1([\P^1/\mu_n],L) =
H^1(\P^1,L)^{\mu_n}$.  Consequently $H^1([\P^1/\mu_n],\O(d-e))
=H^1(\P^1,\O(d-e))^{\mu_n} = 0$, and hence $\tilde{E}$ splits over
$[\P^1/\mu_n]$ as a sum $\O(d) \oplus \O(e)$.  This reduces its
structure group to the maximal torus $\tilde{T}$ of $GL_2$,
determining a $\mu_n$-invariant section of $\tilde{E}/\tilde{T} =
E/T$, which in turn reduces $E$ to $T$ as a bundle over $[\P^1/\mu_n]$.

{\bf Steps B, C, and D} then proceed exactly as in the proof of
Theorem \re{line-exist}.  \fp

\bs{Corollary}
\Corollary{[\P^1/\mu_n]}
\es

\pf.  Similar to that of Corollary \re{line-corollary}.  \fp

\bs{Lemma}
\label{splitting-lemma}
Suppose that a diagonalizable group scheme $S$ acts on a scheme $X$
satisfying $\Ga(\O_X^\times) = (\Gm)_k$, and suppose that $\Pic X$ is
\'etale, generated by a finite number of line bundles on which the
$S$-action is linearized.  Then the short exact sequence
$$1 \lrow \hat{S} \lrow \Pic [X/S] \to \Pic X \to 1,$$
where $\hat{S} = \Hom(S,\Gm)$, 
is split by sending a line bundle $L$ over $X$ to its unique lifting
to a line bundle over $[X/S]$ such that $S$ acts trivially on the
fiber over $x$.  \es

\pf. The last hypothesis guarantees that the natural homomorphism
$\Pic [X/S] \to \Pic X$ is surjective, and the first guarantees that
its kernel, consisting of the group of linearizations on $\O_X$, is
isomorphic to $\hat{S}$.  The \'etale
hypothesis then guarantees that the splitting described is a morphism.
\fp

\bs{Theorem (Uniqueness for {\boldmath $\mu_n$}-equivariant lines)}
\label{mu-line-unique}
\Unique{[\P^1/\mu_n]}
\es

\pf.  Similar to that of Theorem \re{line-unique}, but with the
following additions.

For a character $\chi \in \hat{\mu}_n$ and an integer $d \in \Z$,
denote by $\O(\chi,d)$ the lifting of $\O(d)$ to $[\P^1/\mu_n]$ such
that $\mu_n$ acts on the fiber over $[1,0]$ via $\chi$.  By Lemma
\re{splitting-lemma}, every line bundle over $[\P^1/\mu_n]$ is
isomorphic to exactly one of these.

In the case $G = GL_n$, therefore, one has to prove that if
$\bigoplus\O(\chi_i,a_i) = \bigoplus \O(\psi_j,b_j)$, then the
$(\chi_i,a_i)$ and the $(\psi_j,b_j)$ are the same up to permutation.
This is accomplished by arguing as for Theorem \re{line-unique} that
there exist summands $\O(\psi_j,b_j)$ of degree $a_n$, then using
Schur's lemma to find one such summand where $\psi_j = \chi_n$.

For general $G$, assume first that $k$ has characteristic zero.  There
is a $W$-equivariant isomorphism $P \otimes \La \cong
\Hom(\hat{P},T)$, so regard $\la, \la' \in P \otimes \La$ as
homomorphisms $\hat{P} \to T$.  Since $P = \Pic [\P^1/H] \cong \Z
\times \Z/n\Z$, its Cartier dual satisfies $\hat{P} \cong \Gm \times \mu_n$ and
hence is a smooth group scheme having a $k$-rational point $t_0$ that
generates a dense subgroup.  Then follow the proof of Theorem
\re{line-unique}.

When $k$ has characteristic $p>0$, there is the complication that
$\hat{P} \cong \Gm \times \mu_n$ may not be reduced, so homomorphisms
are not determined by their values on geometric points.  Suppose,
however, that $\la, \la' \in P \otimes \La$ determine isomorphic
$G$-bundles on $[\P^1/\mu_n]$.  If $G \to GL_N$ is any representation
over $k$, then the associated $GL_N$-bundles are also isomorphic and
therefore correspond to the same element of $(P \otimes \La_N)/W_N$,
where $\La_N$ and $W_N$ are the cocharacter lattice and Weyl group of
$GL_N$.

Now recall that every pair $(G,T)$ consisting of a reductive group and
maximal torus over $k$ is the base change of a Chevalley group $G_\Z$
over $\Z$ with maximal torus $T_\Z$, and that the cocharacter lattice
$\La_\Z = \Hom (\Gm, T_\Z)$ is constant over $\Spec \Z$ as a group
with $W$-action.  For any dominant weight $\be$, let $L_\be$ be the
corresponding line bundle over $G_\Z/B_\Z$, and let $V_\be = H^0(G_\Z
/ B_\Z, L_\be)$.  Then $V_\be$ is a finitely generated, free abelian
group acted on by $G_\Z$, and for any field $K$, the Borel-Weil
representation satisfies $H^0(G_K/B_K,L_\be) = V_\be \otimes K$.
Choosing an integer basis for $V_\be$ consisting of weight vectors
determines a homomorphism $G_\Z \to (GL_N)_\Z$ taking $T_\Z$ to the
diagonal matrices.  The induced map $\La \to \La_N$ is therefore
constant, hence independent of $K$.

For $K = \C$, all irreducible representations appear in this fashion.
For $K = k$, our field of characteristic $p$, the Borel-Weil
representation may be reducible; nevertheless, the elements of $(P
\otimes \La_N)/W_N$ corresponding to $\la$ and $\la'$ must be the
same.  The map $(P \otimes \La)/W \to (P \otimes \La_N)/W_N$ is the
unique one making the diagram
$$\begin{array}{ccc}
P \otimes \La & \lrow & P \otimes \La_N \\
\down && \down \\
(P \otimes \La)/W & \lrow & (P \otimes \La_N)/W_N
\end{array}$$
commute, and the top row is independent of $K$, so the bottom row is
as well.  Hence for $K = \C$ as well as $K = k$, the elements of $(P
\otimes \La_N)/W_N$ corresponding to $\la$ and $\la'$ are equal.  By
the characteristic zero case, the $G$-bundles over $\C$ corresponding
to $\la$ and $\la'$ must be isomorphic, so there exists $w \in W$ such
that $\la' = w \la$, as desired.  \fp

\br{Remark.}  
\label{no-harder-narasimhan}
There is in general no unique Harder-Narasimhan
reduction in the sense of Theorem \re{harder-narasimhan} for bundles
on $\mu_n$-equivariant lines.  Even on $[\P^1/\mu_2]$, if one lets $E$
be $\O \oplus \O \to \P^1$ where $\mu_2$ acts with weights $1$ and
$-1$ on the two summands, then either summand gives a reduction to $B
\subset GL_2$ satisfying all the conditions of Theorem
\re{harder-narasimhan}.  That is why the following proof uses the
Harder-Narasimhan reduction on $\P^1$ itself rather than on
$[\P^1/\mu_n]$.  \er

\bs{Theorem (Connectedness for {\boldmath $\mu_n$}-equivariant lines)}
\label{mu-line-connect}
Suppose the derived subgroup of $G$ is simply connected.  Then the
\Connect{[\P^1/\mu_n]}
\es

\pf.  Pull back $E$ from $[\P^1/\mu_n]$ to $\P^1$ and let $E_P$ be the
Harder-Narasimhan reduction of Theorem \re{harder-narasimhan}.  An
automorphism of $E$ over $[\P^1/\mu_n]$ is the same as an automorphism
over $\P^1$ commuting with the $\mu_n$-action.  Since the
Harder-Narasimhan reduction is unique and rigid, it must be preserved
both by the automorphism and by the $\mu_n$-action.  Hence $\Ad E_P$
descends to $[\P^1/\mu_n]$, and there $\Ga(\Ad E) = \Ga(\Ad E_P)$.

As in the proof of Theorem \re{line-connected}, decompose $P = L
\ltimes U$, a semidirect product of a Levi factor and a maximal
unipotent, and let $\Ad E_P = Q \ltimes R$ be the corresponding
decomposition.  Since $\mu_n$ acts on $E$ via a reduction to $T$ by
Theorem \re{mu-line-exist}, it preserves this decomposition, as well
as the splitting of $R$ into a product of line bundles.  So $Q$ and
$R$ (and its splitting) descend to $[\P^1/\mu_n]$.

Consider first $\Ga(Q)$.  The action of $\mu_n$ on the trivial bundle
$Q$ over $\P^1$ must be by conjugation by some homomorphism $\mu_n \to
L$.  Hence $\Ga(Q) = Z_L(\mu_n)$, the centralizer of the image of this
homomorphism.  As $L$ is a torus centralizer in a reductive group $G$
with simply connected derived group, it too has simply connected
derived group: see Steinberg \cite[2.17]{steinberg3}.  But the
centralizer of the image of $\mu_n$ in such a group is smooth \cite[XI
5.2]{sga3} and connected.  If the image of $\mu_n$ is \'etale, as in
characteristic zero, the connectedness is a well-known result of
Steinberg \cite[8.1]{steinberg2}.  The necessary argument in the
general case was kindly communicated to us by Brian Conrad and is set
down in Theorem \re{conrad-theorem} of the Appendix.

As for $\Ga(R)$, the argument in the proof of Theorem
\re{line-connected} goes through without change, showing that $\Ga(R)$
is smooth, affine, and connected.

The kernel of the evaluation map is treated similarly, regardless of
whether the rational points are fixed by $\mu_n$.  \fp

\br{Remark.} In the preceding proof, and throughout the paper,
centralizers are taken in the scheme-theoretic sense.  This should cause
no confusion, in light of the aforementioned result \cite[XI
  5.2]{sga3} implying that the centralizer, in a linear algebraic
group, of a diagonalizable sub-group scheme is a smooth sub-group
scheme, that is to say, a linear algebraic subgroup.

\br{Counterexample} where $G$ is not simply connected.  Let $\mu_2 \to
PGL_2$ be a nontrivial homomorphism extending to $\Gm \to PGL_2$.
This determines a $PGL_2$-bundle over $[\P^1/\mu_2]$ whose pullback to
$\P^1$ is trivial.  By Lemma \re{mu-line-lemma} it is rationally
trivial.  But its automorphism group is the centralizer
$Z_{PGL_2}(\mu_2)$, which is disconnected.

\bit{Variation \thesubsection: a football}
\newcounter{football}
\setcounter{football}{\thesubsection}
Now let $a,b$ be positive integers.  Let $\P^1_{a,b}$ be the smooth
stack, tame in the sense of Abramovich-Olsson-Vistoli \cite{aov},
whose coarse moduli space is $\P^1$, but with isotropy $\mu_a$ over
$p_+ = [1,0]$, isotropy $\mu_b$ over $[0,1]= p_-$, and trivial
isotropy elsewhere.  We call such a stack a {\it football} (referring
to an American football, not a soccer ball).  It can be constructed by
performing the root construction of Cadman \cite{cadman} on $\P^1$,
first to order $a$ at $p_+$, then to order $b$ at $p_-$.  Note that
the global quotient $[\P^1/\mu_{ab}] = \P^1_{ab,ab}$ admits a morphism
to $\P^1_{a,b}$ inducing an isomorphism on coarse moduli spaces.

\bs{Lemma} 
\label{descent-lemma}
The morphism above induces an equivalence between the category of all
vector bundles over $\P^1_{a,b}$ and the category of those vector
bundles over $\P^1_{ab,ab}$ on which the subgroups $\mu_b$ and $\mu_a$
of the isotropy act trivially at $p_+$ and $p_-$, respectively.  
\es

\pf.  
Let $\phi: [\A^1/\mu_a] \to \A^1/\mu_a$ be the coarse
moduli map.  Alper proves \cite[4.5]{alper} that for any vector bundle
$E$ over the target, the natural adjunction map $E \to \phi_* \phi^*
E$ is an isomorphism.  He also proves \cite[10.3]{alper} that, for any
vector bundle $F$ over the source on which the isotropy group at the
origin acts trivially, the natural adjunction map $\phi^* \phi_* F \to
F$ is an isomorphism.  

The standard action of $\mu_{ab}$ on $\A^1$ induces an action of
$\mu_b \cong \mu_{ab}/\mu_a$ on $\A^1/\mu_a$.  If a vector
bundle $E$ or $F$ is linearized for this action, then so are its
pullbacks and pushforwards, and the adjunction maps are compatible
with these linearizations.  Hence pullback and pushforward by the
corresponding map $[\A^1/\mu_{ab}] \to [\A^1/\mu_b]$ induces an
equivalence between the category of all vector bundles over
$[\A^1/\mu_b]$ and the category of those vector bundles over
$[\A^1/\mu_{ab}]$ on which the subgroup $\mu_a$ of the isotropy at the
origin acts trivially.

Since a similar statement holds for $[\A^1/\mu_{ab}] \to
[\A^1/\mu_a]$, and the assertion is local on the base, the result
follows.  \fp

\bs{Proposition}
\label{football-picard}
The Picard group $\Pic \P^1_{a,b}$ is the free abelian group generated by
$\O(p_+)$ and $\O(p_-)$, modulo the subgroup generated by $\O(ap_+-bp_-)$.
\es

\pf.  Straightforward using Lemma \re{descent-lemma} for line bundles
and then applying Lemma \re{splitting-lemma} to $[\P^1/\mu_{ab}]$.
\fp

Define the {\it degree} to be the map $\deg: \Pic \P^1_{a,b}
\to \Q$ given by $\deg \O(ip_+ + jp_-) = i/a + j/b$.

\bs{Proposition} 
\label{football-cohomology}
We have $h^0(\P^1_{a,b},\O(ip_+ + jp_-)) = \max (\lfloor i/a
\rfloor + \lfloor j/b \rfloor + 1,0)$.  In particular, $\deg <
0 $ implies $h^0 = 0$.  There exists a section
nonvanishing at $p_+$ (resp.\ $p_-$) if and only if $i|a$
(resp.\ $j|b$).  The canonical bundle is $\kappa_{\P^1_{a,b}}
\cong \O(p_+ + p_-)$, so by Serre duality $h^1(\P^1_{a,b},\O(ip_+
+ jp_-)) = h^0(\P^1_{a,b},\O((1-i)p_+ + (1-j)p_-))$.  In
particular, $\deg \geq 0 $ implies $h^1 = 0$.
\es 

\pf.  Again, straightforward from Lemma \re{descent-lemma}.  \fp

A principal $G$-bundle over $\P^1_{a,b}$ is said to be {\it rationally
  trivial\/} if its restriction to the generic point $\Spec k(t)$ of
$\P^1_{a,b}$ is trivial.

\bs{Theorem (Existence for footballs)} 
\label{football-exist}
\Exist{\P^1_{a,b}}
\es

\pf.  
{\bf Step A}: {\it The case $G=PGL_2$.}
Arguing exactly as in the proof of Theorem \re{line-exist} reduces the $PGL_2$
case to the $GL_2$ case.

For $G= GL_2$, by Lemma \re{descent-lemma} it suffices to consider the
case of $\P^1_{n,n} = [\P^1/\mu_n]$ for $n = ab$; but this is covered
by Theorem \re{mu-line-exist}.

{\bf Steps B, C, and D} then proceed exactly as in the proof of
Theorem \re{line-exist}, using Lemma \re{football-cohomology} for the
vanishing of $H^1$ in Step D.  \fp

\br{Counterexample} where there are three orbifold points on the line.
Let $X = \P^1_{2,2,2}$, that is, the stack obtained from $\P^1$ by
performing the root construction to order 2 at each of 3 rational
points.  As shown by Borne \cite[3.13]{borne}, there is an equivalence
of categories between vector bundles on $X$ and parabolic bundles on
$\P^1$ with 3 marked points having half-integral weights.  So let $E =
\O \otimes k^2$ be the trivial bundle on $\P^1$, endowed with a
parabolic structure at the 3 marked points consisting of 3 distinct flags
in $k^2$, each with weights $0$ and $1/2$.  Clearly this does not
split as a sum of parabolic line bundles, so the corresponding vector
bundle on $X$ does not split either. \er

\bs{Corollary}
\Corollary{\P^1_{a,b}}
\es

\pf.  Similar to that of Corollary \re{line-corollary}. \fp

\bs{Theorem (Uniqueness for footballs)}
\label{football-unique}
\Unique{\P^1_{a,b}}
\es

\pf.  The case $G = GL_n$ is an immediate consequence of Theorem
\re{mu-line-unique} and Lemma \re{descent-lemma}.  The general case
follows from this case exactly as in the proof of Theorem
\re{mu-line-unique}.  \fp

\bs{Theorem (Connectedness for footballs)}
\label{football-connect}
Suppose the derived subgroup of $G$ is simply connected.  Then the
\Connect{\P^1_{a,b}}
\es

\pf.  
By Theorem \re{football-exist}, $E$ reduces to $T$.  Hence the
splitting of ${\mathfrak g}$ as a $T$-representation into Cartan
subalgebra and root spaces induces a splitting
$$\ad E \cong \O^r \oplus \bigoplus_{\al \in \Phi} L_\al.$$

As seen in the proof of Theorem \re{harder-narasimhan}, a reduction
$E_P$ of $E$ to a parabolic $P$ is determined by a bundle $\ad E_P$ of
subalgebras of $\ad E$ conjugate to ${\mathfrak p}$.  As was done
there, let $\ad E_P$ be the bundle consisting of line bundles of
degree $\geq 0$ in the splitting above.  In light of Remark
\re{no-harder-narasimhan}, this may not be the smallest such rigid
reduction; nevertheless, it follows from Proposition
\re{football-cohomology} that any automorphism of $\ad E$ preserves
$\ad E_P$, and hence that $\Ga (\Ad E) = \Ga (\Ad E_P)$, since the
normalizer of ${\mathfrak p} \subset {\mathfrak g}$ in $G$ is $P$.  Since $T$
preserves the semidirect product decomposition $P = L \ltimes U$ of
the corresponding parabolic, we have a decomposition $\Ad E_P = Q
\ltimes R$.  Hence as a scheme $\Ga (\Ad E_P) = \Ga(Q) \times \Ga(R)$.

Consider first $\Ga(Q)$.  By Proposition \re{football-picard}, $\Pic
\P^1_{a,b} \cong \Z \times \Z/d\Z$, where $d = (a,b)$ is the greatest
common divisor, and the factor $\Z/d\Z$ corresponds to the bundles of
degree zero.  Pullback by the $\mu_d$-cover $\pi:\P^1_{a/d,b/d} \to
\P^1_{a,b}$ kills this $\Z/d\Z$ factor.  Hence $\pi^* {\mathfrak q}$
is trivial. Since $\pi^* Q$ is an adjoint bundle for the connected
group $L$, this implies that $\pi^* Q$ is also trivial.  Hence
$\Ga(\pi^* Q) \cong L$, which, being a torus centralizer in $G$, has
simply connected derived group \cite[2.17]{steinberg3}.  Then $\Ga(Q)$
is the invariant part under the $\mu_d$-action, namely the centralizer
$Z_L(\mu_d)$.  This is smooth \cite[XI 5.2]{sga3} and, of course,
affine.  By Theorem \re{conrad-theorem} from the Appendix, it is
connected.

As for $\Ga(R)$, observe that since $E$ reduces to $T$, and the
$T$-action on $U$ preserves the root subgroups that directly span it,
the bundle $R$ splits as a product of line bundles, as before.  Hence the
argument in the proof of Theorem \re{line-connected} goes through
without change.  

The kernel of the evaluation map is treated similarly, regardless of
whether any of the rational points coincide with $p_+$ or $p_-$.  \fp

\bit{Variation \thesubsection: a gerbe over a football}
We now extend the results of Variation \thefootball\ to bundles over a
more general base space: a rationally trivial $\mu_n$-gerbe $J$ over a
football, whose isotropy groups at $p_+$ and $p_-$ are $\mu_{na}$ and
$\mu_{nb}$, respectively.  Here, by a {\it $\mu_n$-gerbe} over the
football, we mean a stack over the football, locally isomorphic to a
product with $B \mu_n$ in, say, the fppf topology on the football.
Before we proceed, here are two counterexamples indicating why the
hypothesis on the isotropy is necessary.

\br{Counterexample} where the isotropy hypothesis does not hold.   Over
any field with two distinct roots $\pm i$ of $-1$, let $X =
[\P^1/\mu_2] \times B\mu_2$.  This is a trivial gerbe over
$[\P^1/\mu_2] = \P^1_{2,2}$ with structure group $\mu_2$.  The
quaternion group $\{ \pm 1, \pm i, \pm j, \pm k \}$, acting on a
2-dimensional vector space with basis $\{1,j\}$, defines an
irreducible 2-dimensional projective representation of $\mu_2^2$,
hence a $PGL_2$-bundle over $X$ which does not reduce to a maximal
torus. \er

\br{Counterexample} where the band is nontrivial.  Over any field with
two distinct roots $\pm i$ of $-1$, let $Y = [\P^1/D_4]$ where the
dihedral group $D_4$ acts via its quotient $\mu_2 = \Z/2\Z$.  This is
a $\mu_4$-gerbe over the football $[\P^1/\mu_2]$ with nontrivial band.
Any irreducible representation of $D_4$ of dimension $2$ induces a
vector bundle over $Y$, which does not split into a sum of line
bundles; this may be seen by restricting it to one of the fixed
points. \er

On the other hand, our hypothesis on the isotropy guarantees better
behavior.

\bs{Proposition} 
\label{gerbe-trivial-band}
Any $\mu_n$-gerbe $J$ over a football $\P^1_{a,b}$
with abelian isotropy at $p_+$ and $p_-$ must have trivial band. \es

\pf.  Since $\Aut \mu_n = \Aut \hat{\mu}_n = \Aut \Z/n\Z =
(\Z/n\Z)^\times$, the isomorphism class of the band lies in
$H^1(\P^1_{a,b},(\Z/n\Z)^\times)$.  By hypothesis, the isotropy groups
at $p_\pm$ are central extensions of $\mu_a$ and $\mu_b$ by $\mu_n$.
Hence the gerbes obtained by restricting $J$ to $B\mu_a$ and $B\mu_b$
have trivial band.  It therefore suffices to show that the restriction map
$H^1(\P^1_{a,b},(\Z/n\Z)^\times) \to H^1(B\mu_a,(\Z/n\Z)^\times) \times
H^1(B\mu_b, (\Z/n\Z)^\times)$ is injective.  Since $(\Z/n\Z)^\times$ is
a finite abelian group, it suffices to prove the following lemma.  \fp

\bs{Lemma}
\label{gerbe-lemma}
The restriction
 $H^1(\P^1_{a,b},
\Z/q\Z) \to H^1(B\mu_a, \Z/q\Z) \times H^1(B\mu_b, \Z/q\Z)$ is
injective for any prime power $q$.
\es

\pf. When $q$ is not a power of the characteristic, this is easy, for then 
$$0 \lrow \Z/q\Z \lrow \Gm \lrow \Gm \lrow 1.$$
We therefore have a diagram
$$\begin{array}{ccc}
H^1(\P^1_{a,b}, \Z/q\Z) & \lrow & H^1(\P^1_{a,b}, \Gm) \\
\down & & \down   \\
H^1(B\mu_a \cup B\mu_b, \Z/q\Z) & \lrow & H^1(B\mu_a \cup B\mu_b, \Gm)
\end{array} 
$$
whose horizontal arrows are injective.  By Theorem 90 for stacks, the
spaces on the right-hand side are Picard groups.  By Proposition
\re{football-picard} the right-hand map can be identified with the
obvious map $\Z^2/\langle (a,-b) \rangle \to \Z/a\Z \times \Z/b\Z$.
This is not injective, but it becomes injective when restricted to the
kernel of multiplication by $q$, which is the image of the horizontal
map.  The left-hand vertical map is therefore injective.

On the other hand, if $q = \Char k$, we have the Artin-Schreier
sequence 
$$0 \lrow \Z/q\Z \lrow {\mathbb G}_a \stackrel{F-1}{\lrow} {\mathbb
  G}_a \lrow 0.$$ 
From its long exact sequence, together with $H^1(\P^1_{a,b},\O) =0$ from
Proposition \re{football-cohomology}, we deduce that $H^1(\P^1_{a,b},
\Z/q\Z) = 0$.  More generally, if $q$ is a power of $\Char k$, then we use 
an analogue of
the Artin-Schreier sequence, namely
$$0 \lrow \Z/q\Z \lrow W_i \stackrel{F-1}{\lrow} W_i \lrow 0$$
where $W_i$ is the sheaf of Witt vectors \cite{witt}.  Since these also
satisfy
$$0 \lrow \O \lrow W_i \lrow W_{i-1} \lrow 0,$$
the vanishing of $H^1(\P^1_{a,b},\Z/q\Z)$ again follows from that
of $H^1(\P^1_{a,b},\O)$.  \fp

However, the isotropy hypothesis is not enough.  The gerbe $J$ is said
to be {\it rationally trivial} if its restriction to the generic point
$\Spec k(t)$ of $\P^1_{a,b}$ is trivial.  (By Tsen's theorem 
this is automatic if $k$ is algebraically closed.)  Something like
this is necessary, as the following example shows.

\br{Counterexample} where the gerbe is not rationally trivial.  Since
the absolute Galois group of $\R$ is $\Z/2\Z$, we have $H^2(\Spec \R,
\Z/2\Z) = H^2(\Z/2\Z,\Z/2\Z) = \Z/2\Z$.  Hence there is a nontrivial
$\mu_2 = \Z/2\Z$-gerbe over $\Spec \R$.  Indeed, it is the quotient
$[\Spec \C/\mu_4]$, where the generator of $\mu_4$ acts by complex
conjugation.  The quaternions, with the generator of $\mu_4$ acting by
multiplication by $j$, define a $GL_2$-bundle over this stack which is
irreducible and hence does not reduce to $T$.  The same is true of its
pullback under $\P^1_{a,b} \to \Spec \R$.

\bs{Proposition}
\label{gerbe-root}
Any rationally trivial $\mu_n$-gerbe over a football $\P^1_{a,b}$ with
trivial band is the gerbe of $n$th roots of a line bundle.
\es

\pf.  The short exact sequence
$$1 \lrow \mu_n \lrow \Gm \lrow \Gm \lrow 1$$
gives long exact sequences (say in fppf cohomology)
$$ \begin{array}{ccccccc}
H^1(\P^1_{a,b},\Gm) & \lrow & 
H^1(\P^1_{a,b},\Gm) & \lrow & 
H^2(\P^1_{a,b},\mu_n) & \lrow & 
H^2(\P^1_{a,b},\Gm) \\
\down && \down && \down && \down \\
H^1(\eta,\Gm) & \lrow & 
H^1(\eta,\Gm) & \lrow & 
H^2(\eta,\mu_n) & \lrow & 
H^2(\eta,\Gm), 
\end{array} $$
where $\eta$ denotes the generic point.  The images of the middle
horizontal arrows are the root gerbes.  According to Lieblich
\cite[3.1.3.3]{lieblich}, the right-hand vertical arrow is injective. \fp

\bs{Proposition} 
\label{gerbe-cohomology}
Let $J$ be the gerbe of $n$th roots of a line bundle $L$ over
$\P^1_{a,b}$.  The Picard group of $J$ is generated over $\Pic
\P^1_{a,b}$ by the tautological line bundle $\xi$ modulo the single
relation $\xi^n \cong L$.  In particular, there is a well-defined
degree $\Pic J \to \Q$ extending the degree on $\Pic \P^1_{a,b}$.  The
cohomology of a line bundle $M$ is $H^i(J,M) \cong H^i(\P^1_{a,b},N)$
if $M$ is the pullback of a line bundle $N$ over $\P^1_{a,b}$, and $0$
otherwise.  In particular, $\deg \geq 0$ implies $H^1 = 0$, whereas
$\deg \leq 0$ implies $H^0 = 0$ except for a trivial bundle.  \es

\pf.  The first statement is proved by Cadman \cite[3.1.2]{cadman}.
The second follows from applying the Leray sequence to the projection
$J \to \P^1_{a,b}$.  \fp

A $G$-bundle over the rationally trivial gerbe $J$ is said to be {\it
  rationally trivial} if, on the generic point $\Spec k(t) \times
B\mu_n$ of $J$, it is a trivial $G$-bundle with $\mu_n$ acting by a
homomorphism into a split torus of $G_{k(t)}$.  The last condition may seem ad hoc
(and is non-vacuous even if $k$ is algebraically closed), but
something of this sort is necessary to avoid counterexamples similar
to Counterexample \re{mu-line-example}, only with $\mu_2$ acting trivially on
$\P^1$.

\bs{Theorem (Existence for gerbes)}
\label{gerbe-exist}
Let $J$ be a rationally trivial $\mu_n$-gerbe over $\P^1_{a,b}$,
having isotropy $\mu_{na}$ and $\mu_{nb}$ over $p_+$ and $p_-$,
respectively.
\Exist{J}
\es

\pf.  
\nopagebreak

{\bf Step A}: {\it The case $G=PGL_2$.}  Suppose $E \to J$ is a rationally
trivial principal $PGL_2$-bundle.  Let $B\mu_2 \to J' \to J$ be the
gerbe of liftings of $E$ to $SL_2$.  The pullback of $E$ to $J'$
tautologically lifts to an $SL_2$-bundle $E'$.  

But $J'$ is also a gerbe over the same football, whose structure
group, say $H$, is an extension of $\mu_n$ by $\mu_2$.  As $\Aut
\mu_2 = \Aut \hat{\mu}_2 = \Aut \Z/2\Z = 1$, all such extensions are
central.  Since $\Hom(\mu_n,\mu_2) = \Hom(\hat{\mu}_2,\hat{\mu}_n)$ is
either $\Z/2\Z$ or trivial, the commutator map $\mu_n \to
\Hom(\mu_n,\mu_2)$ factors through the group of components of $\mu_n$,
which is cyclic.  But since every central extension of a cyclic group
is abelian, this implies that $H$ is abelian, and hence is either
$\mu_{2n}$ or $\mu_2 \times \mu_n$.  By the same token, the isotropy
groups of $J'$ over $p_\pm$ are also abelian.  If $H \cong \mu_{2n}$,
then $J'$ has trivial band by Proposition \re{gerbe-trivial-band}.

If $H \cong \mu_2 \times \mu_n$, then the band of $J'$ is a bundle
over $\P^1_{a,b}$ of groups isomorphic to $H$, with structure group
$\Aut H$.  In fact, as $J'$ is a gerbe over a gerbe, this structure
group reduces to the subgroup preserving $\mu_2$.  The band surjects
onto the band of $J$, which is the trivial $\mu_n$-bundle, and its
kernel is a bundle of groups isomorphic to $\mu_2$, necessarily
trivial as well since $\Aut \mu_2 = 1$.  Hence the structure
group of the band further reduces to $\Hom (\mu_n,\mu_2)$, which is
$\Z/2\Z$ or trivial.  Then by Lemma \re{gerbe-lemma}, $J'$ again has
trivial band.

Its structure group $H$ therefore globally acts on $E'$, regarded as a
rank 2 vector bundle.  If there are two distinct characters, this
splits $E'$ globally as a sum of line bundles, reducing it to a
maximal torus $T' \subset GL_2$.

If there is only one character, then let $L$ be a line bundle over
$J'$ on which $H$ acts with the same character.  This exists, as by
Proposition \re{gerbe-root} $J'$ is a root gerbe (or a tensor product
of two root gerbes if $H \cong \mu_2 \times \mu_n$), and root gerbes
admit tautological line bundles on which $\mu_n$ acts by a generator
of $\hat{\mu}_n$.  Then the rank 2 vector bundle $E' \otimes L^{-1}$
is acted on trivially by $H$.

We want to conclude that $E' \otimes L^{-1}$ descends to $\P^1_{a,b}$.
As in the proof of Lemma \re{descent-lemma}, this follows from an
argument like that of Alper \cite[10.3]{alper}.  Over $\P^1_{a,b}$ as
proved in Theorem \re{football-exist}, it reduces to a maximal torus
$T' \subset GL_2$.  Hence so does $E'$ over $J'$.

Since $E'/T' = E/T$ over $J'$, it follows that $E$ reduces to $T$ over
$J'$.  But since $E$ is pulled back from $J$, the isotropy group
$\mu_2$ acts trivially on $E$ and hence on $E/T$ over $J'$.  By
descent for affine morphisms \cite[4.3.1]{vistoli}, the section of
$E/T$ over $J'$ therefore descends to $J$.  Hence $E$ reduces to $T$
over $J$.

{\bf Step B} is a direct consequence of the rational triviality of $E$
over $J$ as we defined it.  Since the generic isotropy $\mu_n$ maps
into a split torus of $G_{k(t)}$ over the generic point $\eta = \Spec
k(t)$ of $\P^1_{a,b}$, and since all maximal split tori are conjugate
to $T_{k(t)}$ \cite[20.9]{borel}, there is a rational section of $E/B$
preserved by the generic isotropy, that is, a section of $E/B$ over
the generic point $\Spec k(t) \times B\mu_n$ of $J$.  By the valuative
criterion, this extends to a regular section of $E/B$ over all of $J$.

{\bf Steps C and D} then proceed exactly as in the proof of Theorem
\re{line-exist}, using Proposition \re{gerbe-cohomology} for the
vanishing of $H^1$ in Step D.  \fp

\bs{Corollary}
\Corollary{J}
\es

\pf. Similar to that of Corollary \re{line-corollary}.  \fp

\bs{Theorem (Uniqueness for gerbes)}
\Unique{J}
\es

\pf.  Similar to that of Theorem \re{mu-line-unique}.  

In the case $G = GL_n$, one has to prove that if a vector bundle
splits in two ways as a sum of line bundles, then the isomorphism
classes of the line bundles are the same up to permutation.  This is
accomplished as in the proof of Theorem \re{line-unique}, using the
last statement of Proposition \re{gerbe-cohomology} to show that a
line bundle of maximal degree in the first sum maps isomorphically to
a line bundle of maximal degree in the second sum, and then proceeding by
induction.

To reduce the general case to the $GL_n$ case in characteristic zero,
as in the proof of Theorem \re{mu-line-unique} it suffices to show
that the Cartier dual $\hat{P}$ of $\Pic J$ has a dense
cyclic subgroup.  As $\Pic J$ is abelian and finitely generated by
Proposition \re{gerbe-cohomology}, this amounts to showing that its
torsion part is cyclic.  It follows from Proposition
\re{football-picard} that the restriction $\Tor \Pic \P^1_{a,b} \to \Pic B
\mu_a$ is injective.  By Proposition \re{gerbe-cohomology} there is a
diagram
\begin{equation*}
\label{tor-pic-j}
\begin{array}{ccccccccc}
0 & \lrow & \Tor \Pic \P^1_{a,b} & \lrow & \Tor \Pic J & \lrow & \Z/n\Z \\
&& \down && \down && \Big\Vert \\
0 & \lrow & \Pic B\mu_a & \lrow & \Pic B \mu_{na} & \lrow & \Z/n\Z & \lrow & 0,
\end{array}
\refstepcounter{equation}
\leqno {\mbox{\bf (\theequation)}}
\end{equation*}
where we have used the isotropy hypothesis that the restriction of $J$
to $p_+$ is isomorphic to $B\mu_{na}$.  Hence $\Tor \Pic J$ is a
subgroup of a cyclic group, hence cyclic.

The proof in the characteristic $p>0$ case proceeds exactly as in the proof of Theorem \re{mu-line-unique}.  \fp

\bs{Theorem (Connectedness for gerbes)}
Let $J$ be as in Theorem \re{gerbe-exist}, and suppose the derived
subgroup of $G$ is simply connected.  Then the
\Connect{J}
\es

\pf. The proof runs parallel to that of Theorem \re{football-connect},
except for one detail.  In order to show the connectedness of
$\Ga(Q)$, it suffices to provide a reduced $\mu_d$-cover $\pi:X \to J$
killing the torsion in $\Pic J$.  That it be reduced is necessary to 
ensure that $\Ga(\pi^* Q)$ consists only of $L$.

By Proposition \re{gerbe-cohomology}, $\Pic J$
is a finitely generated abelian group of rank 1, so by
(\re{tor-pic-j}), $\Pic J \cong Z \times \Z/d\Z$ for some $d \in \Z$.
Let $L$ be a line bundle over $J$ whose isomorphism class generates
the second factor, choose a nowhere zero section of $L^{\otimes d}$,
and let $X$ be the inverse image of this section under the $d$th power
map $L \to L^{\otimes d}$.

The pullback by $X \to J$ tautologically kills $L$.  On the other
hand, $X$ is easily seen to be reduced, as follows.  By Proposition
\re{gerbe-root}, the part of $J$ over $\P^1_{a,b} \sans p_-$ is
$[\A^1/\mu_{na}]$ for some linear action of $\mu_{na}$ on $\A^1$, and
then by Proposition \re{gerbe-cohomology}, over this open set the line
bundle $L$ is induced by a character of $\mu_{na}$, which by
(\re{tor-pic-j}) has order exactly $d$.  Hence the part of $X$ over
this open set is isomorphic to $[\A^1/\mu_{na/d}]$, which is clearly
reduced.  A similar argument applies to the open set $\P^1_{a,b} \sans
p_+$.  \fp

\bit{Variation \thesubsection: a chain of lines} \newcounter{chain}
\setcounter{chain}{\thesubsection}
Now let $n$ be a positive integer.
Let $C$ be a nodal chain of $n$ projective lines, where $[0,1]$ on
the first line coincides with $[1,0]$ on the second line, and so on.
Denote $[1,0]$ on the first line by $p_+$ and $[0,1]$ on the last line
by $p_-$.  Let $\O(1)$ over $\P^1$ have the standard trivialization at
$[0,1]$ and $[1,0]$ and, given $d = (d_1,\dots,d_n) \in \Z^n$, define
$\O(d)$ to be the line bundle over $C$ obtained by gluing together the
line bundles $\O(d_i)$ on the $i$th line preserving (the tensor powers
of) these trivializations.  This induces an isomorphism $\Z^n \to \Pic
C$.  

Tensoring with $\La$ induces a further isomorphism $\La^n \to
H^1(C,T)$, since $H^1(C,T) \cong \Pic C \otimes \La$ canonically.
For $\la \in \La^n$, denote by $E(\la)$ the $G$-bundle associated to
the image of $\la$ under the latter isomorphism.  Note that $E(\la)$
too acquires a standard trivialization at $p_\pm$ and at each of the
nodes.  Thus one may define homomorphisms $V_\pm: \Aut E(\la) \to G$
evaluating an automorphism at $p_\pm$ in terms of the standard
trivializations.

\bs{Lemma}
\label{chain-lemma}
{\rm (a)} The image of $V_\pm$ is a parabolic subgroup $P_\pm$ of $G$, the Lie
algebra of whose Levi factor is the direct sum of $\mathfrak t$ and
those root spaces ${\mathfrak g}_\al$ for which all $\al \cdot \la_i =
0$.  {\rm (b)} Let $X$ be a nonempty set of rational points of $C$ and
$\Aut(E(\la),X)$ the subgroup of automorphisms trivial over $X$.  Then
$V_\pm(\Aut(E(\la),X))$ is a smooth unipotent subgroup of $P_\pm$ directly
spanned by its root subgroups.  \es

\pf.  Without loss of generality consider $V_+$.  As shown by Brion
\cite[4.2]{brion}, $\Aut E(\la)$ is a group scheme, locally of finite
type, with Lie algebra $H^0(\ad E(\la))$.  Let us first determine the
Lie algebra of the image, which is the image of the Lie algebra.  This
amounts to finding the image of the evaluation $H^0(\ad E(\la)) \to
{\mathfrak g}$ at $p_+$.

Since the Lie algebra of $G$ splits as a ${\mathfrak t}$-representation 
into root spaces
$${\mathfrak g} = {\mathfrak t} \oplus \bigoplus_{\al \in \Phi}
{\mathfrak g}_\al,$$ 
the vector bundle $\ad E(\la)$ splits under the reduction of structure
group to $T$ as 
$$\ad E(\la) = (\O \otimes {\mathfrak t}) \oplus \bigoplus_{\al \in \Phi} L_\al,$$
where $L_\al \cong \O(\al \cdot \la_1, \dots, \al \cdot \la_n)$.

Clearly the line bundle $\O(d_1, \dots, d_n)$ has a section
nonvanishing at $p_+$ if and only if $(0, \dots, 0) \preccurlyeq (d_1, \dots,
d_n)$ in the lexicographic ordering.  Only finitely many line bundles
$\O(\al \cdot \la_1, \dots, \al \cdot \la_n)$ appear in the sum above.
So choose rational numbers $x_1 >> x_2 >> \cdots >> x_n > 0$; then the
lexicographic inequality holds for these bundles if and only if $0
\leq \sum x_i (\al \cdot \la_i) = \al \cdot \sum x_i \la_i$.  Furthermore,
equality holds in the latter inequality if and only if all $\al \cdot
\la_i = 0$.  The image of $dV_+$ is therefore
\begin{equation*}
\label{evaluation-image}
{\mathfrak t} \oplus \bigoplus_{\al \cdot (\sum x_i \la_i) \geq 0}
{\mathfrak g}_\al.
\refstepcounter{equation}
\leqno {\mbox{\bf (\theequation)}}
\end{equation*}  
This is a parabolic subalgebra ${\mathfrak p}_+$, and its Levi
factor is 
$${\mathfrak t} \oplus \bigoplus_{\al \cdot (\sum x_i \la_i) = 0}
{\mathfrak g}_\al = 
{\mathfrak t} \oplus \bigoplus_{\al \cdot \la_i = 0}          
{\mathfrak g}_\al.$$

In characteristic zero, where all group schemes over a field are
smooth, this immediately implies that the image of $V_+$ contains the
parabolic subgroup $P_+$ whose Lie algebra is ${\mathfrak p}_+$.  

In general, a little more care is needed, since {\it a priori} the
image might not be smooth.  So note instead that, as the Levi factor
consists of those roots $\al$ where all $\al \cdot \la_i = 0$, the
semidirect product decomposition $P = L \ltimes U$ determines a
decomposition $\ad E_P \cong Q \ltimes R$ as in Theorem
\re{line-connected} where $Q$ is a trivial $L$-bundle over $C$ and $R$
is a unipotent bundle.  The evaluation $\Ga(Q) \to L$ is then clearly
surjective; to see that the evaluation $\Ga(R) \to U$ is surjective as
well, observe that the splitting of $U$ (as a scheme) into a product
of root subgroups $U_\al$ determines a similar splitting of $R$, say $R =
\prod_\al R_\al$, and also that $\Ga(R_\al) \to U_\al$ is
surjective as mentioned before. 

Hence, in any characteristic, the image of $V_+$ is a group subscheme
of $G$ containing $P_+$, whose Lie algebra coincides with that of
$P_+$.  It is therefore smooth with identity component $P_+$; but it
must then be connected \cite[11.16]{borel}, hence coincides with
$P_+$.  This proves (a).

As for (b), by cutting the chain at the nodes in $X$ we may assume
without loss of generality that $X$ contains no nodes, hence defines a
Cartier divisor.  The statement then follows by a similar argument,
splitting the vector bundle $\ad(E(\la))(-X)$ as before, then
observing that the intersection of $\Aut(E(\la),X)$ with $Q$ is
trivial, while its intersection with $R$ is again a product of line
bundles corresponding to the roots.  \fp

Say that a $G$-bundle over $C$ is {\it rationally trivial} if it is
trivial over the generic point of each line in $C$, and let
$\bar{H}^1(C,G)$ denote the set of isomorphism classes of rationally
trivial $G$-bundles, as before.

\bs{Theorem (Existence for chains)}
\label{chain-exist}
\Exist{C}
\es

For $k = {\mathbb C}$, this has been proved by Teodorescu \cite{teodorescu}.

\smallskip

\pf\/ by induction on $n$, the case $n=1$ being Theorem
\re{line-exist}.

For $n>0$, we will prove that the structure group may be reduced to $N(T)$
in such a way that it reduces further to $T$ on each irreducible
component of $C$.  This suffices, as the associated principal $N(T)/T$-bundle
is then trivial on each component and hence clearly has a section.

Let $E$ be a rationally trivial principal $G$-bundle over $C$.
Express $C$ as a union of two strictly shorter chains $C_+$ and $C_-$
intersecting in a single point $x$.  By the induction hypothesis, the
restrictions of $E$ to $C_\pm$ reduce to $T$ and hence are isomorphic
to bundles of the form $E(\la_\pm)$ described earlier.  Therefore $E$
may be obtained from $E(\la_\pm)$, equipped with their standard
trivializations at $x$, by using some gluing parameter $g\in G$ to
identify their fibers.  It suffices to show that this parameter may be
moved into $N(T)$ by acting by automorphisms of $E(\la_\pm)$.  Since
the Bruhat decomposition says $G= BWB$, or rather $G = BN(T)B$, it now
suffices to observe that the images of the evaluation homomorphisms
$V_\pm: \Aut E(\la_\pm) \to G$ both contain a Borel containing $T$, as
an immediate consequence of Lemma \re{chain-lemma}.  \fp

\bs{Corollary}
This descends to a natural surjection $\La^n/W \to H^1(C,G)$.
\es

\pf.  Since $\Z^n \cong \Pic C = H^1(C,\Gm)$, tensoring by $\Lambda$
yields $\La^n \cong H^1(C,T)$ as $W$-modules.  If $\la' = w \la
\in \La^n$, then the $T$-bundles corresponding to $\la'$ and $\la$ are
related by extension of structure group by $w:T \to T$.  Since this
extends to an inner automorphism of $G$, the associated $G$-bundles
$E(\la')$ and $E(\la)$ are isomorphic. \fp

\br{Counterexample} where the base curve is not a chain but still has
arithmetic genus zero.  Let $X$ be a nodal curve consisting of 4
projective lines configured like the letter {\sf E}\@.  Let $V$ be a vector
bundle over $X$ constructed by gluing 3 copies of $\O \oplus \O(1)$
(over each of the horizontal lines) to a trivial bundle $\O \otimes
k^2$ (over the vertical line) so that the fibers of $\O(1)$ over the
nodes are glued to 3 distinct lines in $k^2$.  Then $V$ does not split
as a sum of line bundles.  For the only splittings of $\O \otimes k^2$
arise from splittings of $k^2$, whereas any splitting of $\O \oplus
\O(1)$ must include the given $\O(1)$ summand by the last statement of
Theorem \re{harder-narasimhan}.

For another such example, let $Y$ be the curve consisting of 3
projective lines meeting pairwise nodally in a single point with
3-dimensional Zariski tangent space, like the 3 axes in $\A^3$.  Let
$\pi: X \to Y$ be the obvious morphism.  Then $\pi_*V$ is a vector
bundle over $Y$, which may be constructed analogously to $V$, and
which again does not split, for the same reasons.  \er

\bs{Theorem (Uniqueness for chains)}
\label{chain-unique}
\Unique{C}
\es

\pf\/ by induction on $n$, the case
$n=1$ being Theorem \re{line-unique}.

Let $\la$ and $\la' \in \La^n$ be $n$-tuples determining isomorphic
$G$-bundles $E(\la) \cong E(\la')$ over $C$.  It suffices to find $w \in
W$ such that $\la' = w \la$.

Express $C$ as a union of two shorter chains $C_+ \cup C_-$
intersecting at a single node $x$, the first chain having $m$ lines.  Write
$\la = (\la_+,\la_-)$ for $\la_+ \in \La^m$ and $\la_- \in \La^{n-m}$,
and similarly for $\la'$.  Then $E(\la_\pm)$ are the restrictions of
$E(\la)$ to $C_\pm$, equipped with trivializations at $x$,
and similarly for $E(\la')$.  Certainly $E(\la_\pm) \cong E(\la'_\pm)$, so
by the induction hypothesis, there exist $w_\pm \in W$ such that
$\la'_\pm = w_\pm \la_\pm$.

From bundles $E_\pm$ over $C_\pm$ trivialized at $x$, a bundle over
$C$ may be constructed by identifying the trivialized fibers using a
gluing parameter, that is, a group element $g \in G$; denote this
bundle $E_+ \vee_g E_-$.  Then we have isomorphisms
$$E(\la_+) \vee_e E(\la_-) \cong E(\la) \cong E(\la') \cong E(\la'_+)
\vee_e E(\la'_-) = E(w_+ \la_+) \vee_e E(w_- \la_-).$$ 
On the other hand, extension of structure group by the homomorphisms
$w_\pm:T \to T$ induces isomorphisms $E(w_\pm \la_\pm) \cong
E(\la_\pm)$.  These may act nontrivially on the fiber at $x$, so the
gluing parameter must be adjusted to obtain
$$E(w_+ \la_+) \vee_e E(w_- \la_-) \cong E(\la_+) \vee_{w_-^{-1} w_+}
E(\la_-),$$ 
where we abuse notation by using $w_+$ and $w_-$ to denote any of
their representatives in $N(T)$.  Concatenating these two chains of
isomorphisms, we conclude that there exist automorphisms $\phi_+ \in
\Aut E(\la_+)$ and $\phi_- \in \Aut E(\la_-)$ whose evaluations $v_\pm
= V_\pm(\phi_\pm)$ at $x$ satisfy $v_-^{-1} v_+ = w_-^{-1} w_+ \in
N(T)$.

On the other hand, $v_- v_+^{-1} \in P_- P_+$, where $P_\pm$ are the
parabolic subgroups of Lemma \re{chain-lemma}.  If $W_\pm \subset W$
are the Weyl groups of their Levi factors and $\bar{W} = W_-
\backslash W / W_+$, then the Bruhat decomposition is a disjoint union
$G = P_- \bar{W} P_+$ \cite[21.16]{borel}, meaning that the inverse image
under the quotient $\pi: N(T) \to W$ induces a bijection between the
orbits of $W_- \times W_+$ acting on $W$ and those of $P_- \times P_+$
acting on $G$.  Consequently, $P_-P_+ \cap N(T) = \pi^{-1}(W_-)
\pi^{-1}(W_+)$, so there exist $u_\pm \in \pi^{-1}(W_\pm)$ such that
$u_-^{-1} u_+ = v_-^{-1} v_+$.  

Let $w = w_- u_-^{-1} = w_+ u_+^{-1}$.  Since every element of $W_\pm$
fixes $\la_\pm$, we have $w \la_- = w_- u_-^{-1} \la_- = w_- \la_- =
\la'_-$ and $w \la_+ = w_+ u_+^{-1} \la_+ = w_+ \la_+ = \la'_+$, so
$w \la = \la'$, as desired.  \fp

\bs{Theorem (Connectedness for chains)}
\label{chain-connect}
The \Connect{C}
\es

\pf\/ by induction on $n$, the case
$n=1$ being Theorem \re{line-connected}.

Let $X = \{x_1, \dots, x_\ell\}$ be the set of rational points.
Assume the statement for chains shorter than $C$.  Let $E$ be a
rationally trivial $G$-bundle over $C$.  The kernel of the evaluation
map is $\Aut(E,X)$, the group of automorphisms of $E$ trivial over
$X$.  Without loss of generality, we may assume that $X$ contains no
nodes, for if it does, the automorphism group simply splits as a
product of groups of the same form on shorter chains.

Express $C$ as a union of two shorter chains $C_+ \cup C_-$, glued
along a node $x$.  Let $X_+ = X \cap C_+$ and $E_+ = E|_{C_+}$, and
similarly for $X_-$ and $E_-$.  Trivialize $E$ at $x$, so that
automorphisms of that fiber take values in $G$.  Let $G_+$ be the
image of $\Aut(E_+,X_+)$ in $G$ under evaluation at $x$, and similarly
for $G_-$.  Then $\Aut(E,X)$ becomes the fibered product over $G$ of
$\Aut(E_+, X_+)$ and $\Aut(E_-,X_-)$.  As such, it lies in the short
exact sequence
$$1 \lrow \Aut(E_+,X_+ \cup \{ x \}) \times \Aut(E_-,X_- \cup \{ x \})
\lrow \Aut(E,X) \lrow G_+ \cap G_- \lrow 1.$$ 
The left-hand term is smooth, affine, and connected by the induction
hypothesis, so it suffices to show that $G_+ \cap G_-$ is smooth,
affine, and connected.  (That extensions of smooth affine groups are
smooth and affine is clear, as the quotient map is a torsor, locally
trivial in the fppf topology \cite[VIA 3.2]{sga3}, and these
properties are preserved by fppf base change
\cite[1.15]{vistoli}.)


Now apply Lemma \re{chain-lemma}.  If $X$ is empty, it says that
$G_\pm$ are parabolic subgroups, whose intersection is smooth, affine,
and connected \cite[XXVI 4.1.1]{sga3}.  On the other hand, if $X$ is
nonempty, it says that either $G_+$ or $G_-$ is a smooth unipotent
subgroup directly spanned by its root subgroups.  Furthermore, the
other one contains the root subgroup for every root space in its Lie
algebra.  Therefore the intersection is directly spanned by these root
subgroups and hence is smooth, affine, and connected as well.  \fp

\bit{Variation \thesubsection: a torus-equivariant line}
Let $S$ be a split torus acting on $\P^1$ so that $p_+ = [1,0]$ and
$p_- = [0,1]$ are fixed.  This variation is concerned with $G$-bundles
over the stack $[\P^1/S]$, or equivalently, with $S$-equivariant
$G$-bundles over $\P^1$.

A $G$-bundle over $[\P^1/S]$ is said to be {\it rationally trivial\/} if
its pullback to the generic point $\Spec k(t)$ of $\P^1$ is trivial.

\bs{Theorem (Existence for torus-equivariant lines)}
\label{torus-line-exist}
\Exist{[\P^1/S]}
\es

For $G = GL_n$ and $k = {\mathbb C}$, this has been proved by Kumar
\cite{kumar}.

\smallskip

\pf.  {\bf Step A}: {The case $G=PGL_2$.}  This is proved exactly as for
Theorem \re{mu-line-exist}, except that, in the case where the
pullback of $E$ to $\P^1$ is trivial, the result follows immediately,
since the image of the split torus $S$ in $PGL_2$ must be split
\cite[8.4]{borel}, hence either trivial or conjugate to $T$
\cite[20.9]{borel}.

{\bf Step B}: {\it Reduction to $B$.}  Let $E$ be a rationally trivial
$G$-bundle over $[\P^1/S]$.  We distinguish two cases.

{\bf Case 1}: {\it $S$ acts trivially on $\P^1$}.  In this
case, trivializing $E$ over the generic point $\eta = \Spec
k(t)$ gives a homomorphism $\phi: S(t) \to G(t)$ over $k(t)$.
Any split torus acting on a projective variety fixes a rational
point \cite[15.2]{borel}, so $S(t)$ fixes a $k(t)$-rational
point on $G(t)/B(t)$.  Since $B(t)$ is split solvable,
$H^1(k(t),B(t)) = 1$, there is a $k(t)$-rational point $g(t)
\in G(t)$ lying over the fixed point on $G(t)/B(t)$
\cite[15.7]{borel}.  Then $\phi(S(t)) \subset g(t) B(t)
g^{-1}(t)$.  That is to say, $S$ preserves a rational family of
Borel subgroups parametrized by $t$ and defined over an open
subset of $\P^1$.  It therefore preserves a rational section of
$[E/B]$ over $\P^1$.  By the valuative criterion, this extends
to a regular section.

{\bf Case 2}: {\it $S$ acts nontrivially on $\P^1$}.  In this
case, choose a splitting $S \cong S' \times \Gm$ where the torus $S'$ acts
trivially and $\Gm$ acts with weight $n \in \N$.  

The pullback of $E$ to $\P^1$ is rationally trivial; hence by
Theorem \re{line-exist} its restriction to $\A^1$ is trivial.  The
pullback of $E$ to $[\A^1/(S' \times \mu_n)]$ is therefore
induced by a family of homomorphisms $\gamma_t: S' \times \mu_n
\to G$ over $\A^1$.

There is an isomorphism $f:B(S' \times \mu_n) \to [\Gm/S]$;
furthermore, if $i: \Spec k \to \Gm$ is the inclusion of any
$k$-rational point $t$, the diagram
$$\begin{array}{ccc}
\Spec k & \stackrel{i}{\lrow} & \Gm\\
\leftdownarg{P} && \downarg{Q}\\
{[\Spec k/(S'\times \mu_n)]} & \stackrel{f}{\lrow} & [\Gm/S]
\end{array}$$
may be made commutative by choosing the appropriate natural
transformation $f \circ P \Longrightarrow Q \circ i$.  The
restriction of $E$ to $[\Gm/S]$ is therefore a $G$-bundle over
$B(S' \times \mu_n)$ whose pullback to $\Spec k$ is trivial;
indeed, it is induced by the homomorphism $\gamma_t: S' \times \mu_n \to
G$.  

Now $\gamma_0$ extends to a homomorphism $S \to G$.  Its image
therefore lies in a split torus.  Since any family of
homomorphisms from a multiplicative group to a smooth affine
group scheme is \'etale-locally conjugate to a constant family
\cite[2.1.5]{conrad}, this is an open property, so
$\gamma_t$ also has image in a split torus, and hence in a
Borel subgroup.  The bundle $E/B$ therefore has a rational
section which, by the valuative criterion, extends to a regular section.

{\bf Steps C and D} then proceed as in the proof of Theorem
\re{line-exist}.  In Step D, note that since $S$ is a multiplicative
group, its rational representations are completely reducible, so that
taking invariants is exact.  Hence for any line bundle $L$ over
$[\P^1/S]$, $H^1([\P^1/S],L) = H^1(\P^1,L)^S$.  So if
$L$ has positive degree, $H^1([\P^1/\mu_n],L) = 0$ as required. \fp

\bs{Corollary}
\Corollary{[\P^1/S]}
\es

\pf.  Similar to that of Corollary \re{line-corollary}. \fp

\bs{Theorem (Uniqueness for torus-equivariant lines)}
\label{torus-line-unique}
\Unique{[\P^1/S]}
\es

\pf.  Similar to that of Theorem \re{mu-line-unique}, though without
the complication in characteristic $p$.  In the present case the
Cartier dual satisfies $\hat{P} \cong \Gm^{r+1}$, where $r$ is the
rank of $S$; hence, after replacing $k$ by a field extension if necessary,
one may find a $k$-rational point $t_0 \in \hat{P}$ generating a dense
subgroup. \fp

\bs{Theorem (Connectedness for torus-equivariant lines)}
The \Connect{[\P^1/S]}
\es

\pf.  Parallel to that of Theorem \re{line-connected}.  As the
Harder-Narasimhan filtration of Theorem \re{harder-narasimhan} is
unique and rigid, the $S$-action must preserve it.  In the ensuing
decomposition $\Ga(\Ad E) = \Ga(Q) \times \Ga(R)$, the first factor is
the centralizer of a torus in $L$, which is smooth \cite[XI 5.3]{sga3}
and connected \cite[11.12]{borel}, while the second factor is handled
as before.  \fp

\bit{Variation \thesubsection: a torus-equivariant chain of lines}
Let $C$ be a chain of $n$ lines with endpoints $p_\pm$ as in Variation
\thechain.  Let $S$ be a split torus acting on $C$ so that $p_\pm$ are
fixed.  This determines a homomorphism $S \to \Aut(C,p_\pm) \cong
\Gm^n$, not necessarily an immersion.  This variation is concerned with
$G$-bundles over the stack $[C/S]$, or equivalently, with
$S$-equivariant $G$-bundles over $C$.

Given $d \in \Z^n$, a line bundle $\O(d)$ over $C$ of multidegree $d$
was constructed in Variation \thechain, yielding an isomorphism $\Z^n
\to \Pic C$.  Furthermore, given any fixed $k$-rational point $x$ of
the $S$-action on $C$, applying Lemma \re{splitting-lemma} to $[C/S]$
splits the natural exact sequence
\begin{equation*}
\label{torus-chain-sequence}
1 \lrow \hat{S} \lrow \Pic[C/S] \lrow \Pic C \lrow 1,
\refstepcounter{equation}
\leqno {\mbox{\bf (\theequation)}}
\end{equation*}
yielding an isomorphism $\hat{S} \times \Z^n \to \Pic[C/S]$.

Tensoring with $\La = \Hom(\Gm,T)$ induces a further isomorphism
$\Hom(S,T) \times \La^n \to H^1([C/S],T)$.  Once a fixed point $p$ has
been specified, denote by $E(\chi,\la)$ the $G$-bundle associated to
the image of $(\chi,\la) \in \Hom(S,T) \times \La^n$ under the latter
isomorphism.  It inherits from $E(\la)$ a standard trivialization at
$p_\pm$ and at each of the nodes, over each of which $S$ acts on the
trivialized fiber by a homomorphism $S \to T$.  Let $\chi_\pm: S \to
T$ denote these homomorphisms at the points $p_\pm$, and let $S_\pm =
\chi_\pm(S)$.

Any automorphism in $\Aut E(\chi,\la)$ must, by definition, commute
with the $S$-action.  The two homomorphisms $V_\pm: \Aut
E(\chi,\la) \to G$ evaluating an automorphism at the points $p_\pm$
must therefore take values in $Z(S_\pm)$, the centralizers of the
split tori $S_\pm \subset G$.  Note that a root $\al \in \Hom(T,\Gm)$
of $G$ is also a root of $Z(S_\pm)$ if and only if $\al(S_\pm) = 1$,
which is the linear condition that $\al$ lie in the annihilator in
${\mathfrak t}^*$ of ${\mathfrak s}_\pm$.

\bs{Lemma}
\label{equivariant-lemma}
{\rm (a)} The image of $V_\pm$ is a parabolic subgroup $P_\pm$ of $Z(S_\pm)$,
the Lie algebra of whose Levi factor is the direct sum of ${\mathfrak
  t}$ and those root spaces ${\mathfrak g}_\al$ of $Z(S_\pm)$ for
which all $\al \cdot \la_i = 0$.  {\rm (b)} Let $X$ be a nonempty set of rational
points of $C$ and $\Aut(E(\la),X)$ the subgroup of automorphisms
trivial over $X$.  Then $V_\pm(\Aut(E(\la),X))$ is a smooth unipotent subgroup
of $P_\pm$ directly spanned by its root subgroups.  \es

\pf.  The centralizer $Z(S_\pm)$ is smooth \cite[XI 5.3]{sga3},
reductive, and connected \cite[11.12]{borel}.  The proof is therefore
entirely parallel to that of Lemma \re{chain-lemma}.

The only subtlety is in characterizing the sections of the line
bundles $L_\al$ over $[C/S]$.  For such a bundle to have a section
nonvanishing at $p_+$, clearly $S$ must act trivially on $L_\al$ at
$p_+$ (which is the condition that $\al$ be a root of $Z(S_+)$).
Therefore $L_\al$ must be the unique lift of $\O(\al \cdot \la_1,
\dots, \al \cdot \la_n)$ on which $S$ acts trivially over $p_+$.  A
moment's reflection shows that, despite the action of $S$ on $C$, the
unique lift of $\O(d_1, \dots, d_n)$ on which $S$ acts trivially over
$p_+$ has a section over $[C/S]$ nonvanishing at $p_+$ if and only if
$(0, \dots, 0) \preccurlyeq (d_1, \dots, d_n)$ in the lexicographic
ordering.  That is to say, for line bundles on which $S$ acts
trivially at $p_+$, the presence of the $S$-action causes no further
complications.  \fp

\bs{Lemma}
\label{conjugacy-lemma}
Let $G$, $T$, and $S$ be as before.  Any two homomorphisms $S \to T$
which are conjugate by a $k$-rational $g \in G$ are also conjugate by
a $k$-rational $w \in N(T)$.
\es

\pf.  Without loss of generality assume that the homomorphisms are
immersions.  Let their images be $S_\pm \subset T$, with $gS_-g^{-1} =
S_+$.  Then $gTg^{-1} \supset gS_-g^{-1} = S_+$.  Hence the
centralizer $Z(S_+) \supset T \cup gTg^{-1}$, so these two maximal
tori are conjugate in the reductive group $Z(S_+)$ by a $k$-rational
$h \in Z(S_+)$, say $hTh^{-1} = gTg^{-1}$ \cite[20.9]{borel}.  Let $w
= h^{-1}g$.  Then the two homomorphisms are conjugate by $w$, as the
same is true of $g$ while $h^{-1}$ fixes $S_+$; but also $T =
wTw^{-1}$, so $w \in N(T)$.  \fp

A $G$-bundle over $[C/S]$ is said to be {\it rationally trivial} if
its pullback to $C$ is rationally trivial in the sense of Variation
\thechain.

\bs{Theorem (Existence for torus-equivariant chains)}
\label{torus-chain-exist}
\Exist{[C/S]}
\es

\pf\/ by induction on $n$, the case $n=1$ being Theorem \re{torus-line-exist}. 

For $n>0$, let $E$ be a rationally trivial principal $G$-bundle over
$[C/S]$, where $C$ is a chain of $n$ lines.  Express $C$ as a union of
two strictly shorter chains $C_+$ and $C_-$ intersecting in a single
point $x$.  (By the way, even if $S$ acts effectively on $C$, it may
not act effectively on $C_\pm$, which is why we consider this more
general case.)  By the induction hypothesis, the restrictions $E_\pm$
of $E$ to $C_\pm$ may be reduced to $T$ and endowed with their
standard trivializations over the fixed points of $S$.  The actions of
$S$ on the fibers of $E_\pm$ over $x$ are via two homomorphisms
$\chi_\pm: S \to T \subset G$.

These homomorphisms are, of course, conjugate by a $k$-rational
element $g \in G$, namely the relevant gluing parameter.  By Lemma
\re{conjugacy-lemma}, any two such homomorphisms are, in fact,
conjugate by a $k$-rational element $w \in N(T)$.  In light of this,
it clearly suffices to assume that the standard trivializations of
$E_\pm$ over $x$ have the same $S$-action, say by $\chi_+ = \chi_- =
\chi:S \to T$.  For the reduction of $E_-$ to $T$ may be replaced with
its extension of structure group by the homomorphism $w: T \to T$,
which induces an isomorphic $G$-bundle.  The gluing parameter $g \in
G$ with respect to these trivializations commutes with the $S$-action,
hence it belongs to the centralizer $Z(\chi(S))$.

The remainder of the proof exactly follows that of Theorem
\re{chain-exist}, but using the Bruhat decomposition of the reductive
group $Z(\chi(S))$, together with Lemma \re{equivariant-lemma}.  \fp

\bs{Corollary}
Let $P = \Pic [C/S]$; then the surjection above descends to a
natural surjection $(P \otimes \La)/W \to \bar{H}^1([C/S],G)$.  
\es

\pf.  Similar to that of Corollary \re{line-corollary}.  \fp

\bs{Theorem (Uniqueness for torus-equivariant chains)}
\Unique{[C/S]}
\es

\pf\/ by induction on $n$, the case $n=1$ being Theorem \re{torus-line-unique}.

Express $C$ as a union of two shorter chains $C_+ \cup C_-$
intersecting at a single node $x$, the first chain having $m$ lines.
Use the point $x$ to split the exact sequence \re{torus-chain-sequence}
and thereby determine an isomorphism $\Hom(S,T) \times \La^n \to
H^1([C/S],T)$.  For $(\chi,\la) \in \Hom(S,T) \times \La^n$, let
$E(\chi,\la)$ be the corresponding $G$-bundle.

Now suppose that $E(\chi,\la) \cong E(\chi',\la')$.  It suffices to
find $w \in W$ such that $(\chi',\la') = w(\chi,\la)$.  Restricting
both bundles to $x$ shows that $\chi$ and $\chi'$ are conjugate by
some element of $G$, and hence by some element of $W$ according to
Lemma \re{conjugacy-lemma}.  It therefore suffices to assume that
$\chi' = \chi$.

Write $\la = (\la_+,\la_-)$ for $\la_+ \in \La^m$ and $\la_- \in
\La^{n-m}$, and similarly for $\la'$.  Then $E(\chi, \la_\pm)$ are the
restrictions of $E(\chi, \la)$ to $C_\pm$, equipped with given
trivializations at $x$, and similarly for $E(\chi, \la')$.  Certainly
$E(\chi,\la_\pm) \cong E(\chi,\la'_\pm)$, so by the induction
hypothesis, there exist $w_\pm \in W$ such that $(\chi, \la'_\pm) =
w_\pm (\chi, \la_\pm)$.  In particular, $\chi = w_\pm \chi$, so $w_\pm$
belong to $W_{Z(\chi(S))}$, the Weyl group of the centralizer
$Z(\chi(S))$.

The remainder of the proof exactly follows that of Theorem
\re{chain-unique}, but using the Bruhat decomposition of the
centralizer $Z(\chi(S))$, which is smooth \cite[XI 5.3]{sga3},
reductive, and connected \cite[11.12, 13.17]{borel}, together with
Lemma \re{equivariant-lemma}.  \fp

\bs{Theorem (Connectedness for torus-equivariant chains)}
The \Connect{[C/S]}
\es

\pf.  Parallel to that of Theorem \re{chain-connect}, using Lemma
\re{equivariant-lemma}.  When $E$ is trivialized at the node $x$, the
$S$-action on the fiber determines a homomorphism $\phi: S \to G$, and
$G_+$ and $G_-$ are then subgroups of the centralizer $Z_G(\phi(S))$,
which is smooth \cite[XI 5.3]{sga3}, reductive, and
connected \cite[11.12, 13.17]{borel}; it plays the role of $G$ in the
ensuing argument.  \fp

This connectedness theorem, in the case where the torus $S$ is 1- or
2-dimensional, and where the evaluation is at the two endpoints $p_+$
and $p_-$ of the chain, plays a vital part in another paper of the
authors \cite{wonderful}.

\smallskip

There is every reason to expect that the existence, uniqueness, and
connectedness theorems will continue to hold, under suitable
hypotheses, over a common generalization of all the foregoing
variations: say, a torus-equivariant $\mu_n$-gerbe over a chain of
footballs.  You are invited, dear reader, to prove this as an exercise.


\bit{Appendix: connectedness of centralizers in positive characteristic}
We thank Brian Conrad for communicating the following statements and
their proofs.  As in the body of the paper, let $G$ be a connected
reductive group over a field $k$.  (It is elsewhere assumed that $G$
is split, but this is not necessary here.)  Centralizers
are taken in the scheme-theoretic sense, as always.

\bs{Lemma}
\label{conrad-lemma}
For any integer $n \geq 1$, the image of any $k$-homomorphism $\mu_n \to G$ is
contained in a $k$-torus of $G$.  
\es

\pf. Without loss of generality assume that the homomorphism is a
closed immersion, so that $\mu_n \subset G$.

It is equivalent to find a maximal $k$-torus $T$ of $G$ centralizing
$\mu_n$, as such a $T$ is its own centralizer in $G$.  Since the
centralizer $Z_G(\mu_n)$ is smooth \cite[XI 5.3]{sga3}, its maximal
$k$-tori are maximal in $G$ if and only if the same holds true after a
field extension.  Therefore we may assume that $k$ is algebraically
closed.

A theorem of Steinberg \cite[8.1]{steinberg2} then asserts that any
semisimple rational point of $G$ is contained in a torus.  This
implies the lemma whenever the characteristic of $k$ does not divide
$n$, for then $\mu_n$ is \'etale and generated by a semisimple
rational point of $G$.

If the characteristic $p$ does divide $n$, then decompose $\mu_n \cong
E \times \mu_{p^r}$ for some $r$, where $E$ is \'etale.  By
Steinberg's theorem, $E$ is contained in a torus in the
centralizer $Z_G(E)$, which is reductive and contains $\mu_{p^r}$.  It is
therefore in the identity component $H = Z_G(E)^0$, and indeed in its
center.  Replacing $G$ by the quotient $H/E$, we may therefore assume
that $n = p^r$.

Proceed by induction on $r$.  If $r=1$, choose an embedding $G \subset
GL(V)$ as a closed subgroup; since all algebraic representations of
$\mu_p$ are completely reducible, the tangent vector $X$ to $\mu_p$ at
the identity acts semisimply on $k[GL(V)]$, so $X \in {\mathfrak g}$
is a semisimple element \cite[4.3(2), 4.4(4)]{borel}.  It is therefore
tangent to some torus $T \subset G$ \cite[11.8]{borel}.  The $p$-Lie
functor defines a map $\Hom(\mu_p, T) \to
\Hom_{\oper{$p$-Lie}}(\Lie(\mu_p),{\mathfrak t})$, which is bijective
because $\mu_p$ has vanishing relative Frobenius morphism \cite[VIIA
7.2, 7.4]{sga3}.  But the same is true with $T$ replaced by $G$.
Hence the inclusion $\mu_p \subset G$ factors through $T$.

If $r>1$, there is a subgroup $M \subset \mu_{p^r}$ with $M \cong
\mu_p$.  The centralizer $Z_H(M)$ is smooth with
reductive identity component $Z_H(M)^0$ \cite[A.8.12]{cgp}.  Being
central, $M$ lies in a maximal torus of $Z_H(M)^0$, so the inverse
image of some (and hence every) maximal torus in $Z_H(M)^0/M$ is a
maximal torus in $Z_H(M)^0$.  The induction hypothesis then applies to
$\mu_{p^{r-1}} \cong \mu_{p^r}/M \subset Z_H(M)^0/M$, and we are
done. \fp

\bs{Theorem}
\label{conrad-theorem}
Suppose the derived subgroup of $G$ is simply connected.  Then the
cen\-tral\-izer of the image of any homomorphism $\mu_n \to G$ is
connected.
\es

\pf.  Again without loss of generality assume that the homomorphism is
a closed immersion, and also that $k$ is algebraically closed.

Let $M = \mu_n$.  If the characteristic of $k$ does not divide $n$, so
that $M$ is \'etale, then the statement is again a theorem of
Steinberg \cite[2.15]{steinberg3}.

To deduce the general case from this case, first observe that, by
Lemma \re{conrad-lemma}, the subgroup $M \subset G$ lies in a
maximal torus $T$.  

Now recall that every pair $(G,T)$ consisting of a reductive group and
maximal torus over $k$ is the base change of a Chevalley group $G_\Z$
over $\Z$ with maximal torus $T_\Z$.  Every geometric fiber of $G_\Z$
has simply connected derived group, for this property is equivalent to
the quotient of the cocharacter lattice by the coroot lattice being
torsion-free \cite[2.4]{steinberg3}, and these lattices are constant
over $\Z$.

The category of diagonalizable groups over any connected base scheme
is equivalent, under Cartier duality, to the category of abelian
groups.  Both $M$ and $T$ are diagonalizable, so the inclusion $M
\subset T$ is dual to a surjection $\Z^r \to \Z/n\Z$.  Therefore $M$
is the base change of a multiplicative subgroup $M_\Z \subset T_\Z$
dual to the same surjection.

The centralizer $Z_{G_\Z}(M_\Z)$ is then a group subscheme of $G_\Z$,
smooth over $\Z$ \cite[XI 5.3]{sga3}.  Its geometric fibers are the
centralizers of linearly reductive groups and hence have
reductive identity components \cite[A.8.12]{cgp}.

The generic fiber is connected by the theorem of Steinberg mentioned
above.  But the locus of connected fibers, for any smooth affine group
scheme whose fibers have reductive identity components, is always
closed \cite[3.1.12]{conrad}.  Hence the fiber $Z_G(M)$ over $k$
is connected.  \fp

\pagebreak[3]

\end{document}